\documentclass{amsart}

\usepackage{amsmath,amssymb}

\usepackage[all]{xy}

\newtheorem{theorem}{Theorem}[section] 
\newtheorem{lemma}[theorem]{Lemma}     
\newtheorem{lemdfn}[theorem]{Lemma and Definition}
\newtheorem{corollary}[theorem]{Corollary}
\newtheorem{proposition}[theorem]{Proposition}
\newtheorem{definition}{Definition}
\newtheorem{observation}{Observation}

\title[${\rm{PD}}^3$--Pairs]
{Poincar\'e Duality Pairs of Dimension Three} 
 
\author{Beatrice Bleile}

\begin{document}
\maketitle

\begin{abstract}
We extend Hendriks' classification theorem and Turaev's realisation and splitting theorems for ${\rm{PD}}^{3}$--complexes to the relative case of ${\rm{PD}}^{3}$--pairs. The results for ${\rm{PD}}^{3}$--complexes are recovered by restricting the results to the case of ${\rm{PD}}^{3}$--pairs with empty boundary.
Up to oriented homotopy equivalence, ${\rm{PD}}^{3}$--pairs are classified by their fundamental
triple consisting of the fundamental group system, the orientation
character and the image of the fundamental class under the classifying map.
Using the derived module category we provide necessary and
sufficient conditions for a given triple to be realised by a 
${\rm{PD}}^{3}$--pair. The results on classification and realisation yield splitting or
decomposition theorems for ${\rm{PD}}^{3}$--pairs, that is, conditions under
which a given ${\rm{PD}}^{3}$--pair decomposes as interior or boundary
connected sum of two ${\rm{PD}}^{3}$--pairs.

\end{abstract}

\section{Introduction}
A {Poincar{\'e} duality complex of dimension $n$}, or ${\rm{PD}}^n$--complex, is a CW--complex exhibiting $n$--dimensional equivariant Poincar{\'e} duality. We may thus regard ${\rm{PD}}^n$--complexes as {homotopy generalisations} of manifolds. Similarly, {Poincar{\'e} duality pairs of dimension $n$}, or ${\rm{PD}}^n$--pairs, are homotopy generalisations of manifolds with boundary.

A ${\rm{PD}}^n$--pair is a pair of CW--complexes, $(X, \partial X)$, where $X$ is 
connected and $\partial X$ is a ${\rm{PD}}^{n-1}$--complex,
together with an orientation character $\omega \in {\mbox{H}}^1(X; {\mathbb{Z}} / 2 {\mathbb{Z}})$ and a fundamental class $[X, \partial X] \in {\mbox{H}}_3(X, \partial X; {\mathbb{Z}}^{\omega})$, such that 
\[ \cap [X, \partial X]: H^r(X; B) \longrightarrow H_{n-r}(X, \partial X;
B^{\omega})\]
is an isomorphism for every $r \in {\mathbb{Z}}$ and every left module $B$ over the 
integral group ring of the fundamental group of $X$. 

Note that we do not  require $X$ to be finitely dominated. Wall \cite{Wall5} showed that the cellular chain complex $C(X)$ of the universal cover of $X$ is chain homotopy equivalent to a complex of finitely generated projective ${{\Lambda}}$--modules, vanishing except in dimensions $r$ with $0 \leq r \leq n$. Furthermore, $\pi_1(X; \ast)$ is finitely generated and almost finitely presentable, and $X$ is dominated by a finite complex if and only if $\pi_1(X; \ast)$ is finitely presentable.

In Section 2 we review material well-known to experts concerning the definition and properties of the relative twisted cap products but not readily available in the literature, as well as results on the algebraic sums of chain pairs satisfying Poincar\'e duality needed for the discussion of connected sums of Poincar\'e duality pairs. One of the main results is the generalisation of a theorem by Browder \cite{Browder} to the non--simply connected case.

\smallskip

Section 3 is concerned with the homotopy classification of ${\rm{PD}}^3$--pairs. An oriented homotopy equivalence of ${\rm{PD}}^{3}$--pairs, $(X,\partial X)$ 
and $(Y,\partial Y)$, induces an isomorphism 
\[ (\varphi, \{\varphi_{i}\}_{i \in J}):
\{ \kappa_{i}: \pi_1(\partial X_{i}, \ast) \rightarrow \pi_1(X, \ast) \}_{i \in J} \rightarrow 
\{ \rho_{i}: \pi_1(\partial Y_{i}, \ast) \rightarrow \pi_1(Y, \ast) \}_{i \in J}\]
of their $\pi_{1}$--systems such that
\begin{equation}\label{grpsysisoplus1}
\varphi^{\ast}(\omega_{Y}) = \omega_{X} \quad {\rm{and}} \quad
\varphi_{\ast}(c_{X \ast}([X, \partial X])) = c_{Y \ast}([Y, \partial Y]),
\end{equation}
where $c_{X}$ and $c_{Y}$ are classifying maps  and $\{ \partial X_i \}$ and 
$\{ \partial Y_i \}$ are the connected components of $\partial X$ and $\partial Y$ 
respectively. The converse holds for ${\rm{PD}}^3$--pairs with
aspherical boundary components.

Putting $\mu_X := c_{X \ast}([X, \partial X])$, we call $(\{ \kappa_{i} \}_{i \in J}, \omega_X, \mu_X)$ the \emph{fundamental triple} of the ${\rm{PD}}^{3}$--pair $(X,\partial X)$. Two fundamental triples $(\{ \kappa_{i} \}_{i \in J}, \omega, \mu)$ and $(\{ \kappa'_{i} \}_{i \in J}, \omega', \mu')$ are \emph{isomorphic} if there is an isomorphism $(\varphi, \{\varphi_{i}\}_{i \in J}): \{ \kappa_{i} \}_{i \in J} \rightarrow \{ \kappa'_{i} \}_{i \in J}$ with $\varphi^{\ast}(\omega') = \omega$ and $\varphi_{\ast}(\mu) = \mu'$.

\begin{theorem}[{[CLASSIFICATION]}]\label{B}
Two ${\rm{PD}}^{3}$--pairs $(X, \partial X)$ and $(Y, \partial Y)$ with apsherical boundary components are orientedly homotopy equivalent if and only if their fundamental triples are isomorphic.
\end{theorem}
The case $\partial X = \emptyset$ yields Hendriks' Classification Theorem for ${\rm{PD}}^{3}$--complexes.
 
\smallskip
 
Section 4 investigates which triples $(\{\kappa_{i}\}_{i \in J}, \omega, \mu)$ consisting of a $\pi_1$--system $\{\kappa_{i}: G_i \rightarrow G\}_{i \in J}$, $\omega \in {\mbox{H}}^{1}(G; {\mathbb{Z}} / 2 {\mathbb{Z}})$ and $\mu \in {\mbox{H}}_{3}(G, \{G_i\}_{i \in J}; {\mathbb{Z}}^{\omega})$ are realised by ${\rm{PD}}^3$--pairs. We introduce the derived module category \cite{Johnson} which is needed for the formulation of the realisation condition. Given a finitely presentable group $G$ and $\omega \in {\mbox{H}}^{1}(G;{\mathbb{Z}} / 2 {\mathbb{Z}})$, Turaev defined a homomorphism
\[ \nu:  {\mbox{H}}_{3}(G; {\mathbb{Z}}^{\omega}) \longrightarrow [F, I],\]
where $F$ is some ${\mathbb{Z}}[G]$--module, $I$ is the augmentation ideal and $[A, B]$
denotes the group of homotopy classes of ${\mathbb{Z}}[G]$--morphisms from $A$ to $B$. If 
\[ (\{\kappa_{i}: G_i \rightarrow G\}_{i \in J}, \omega, \mu)\]
is the fundamental triple of a ${\rm{PD}}^3$--pair, then $\nu(\mu)$ is a class of homotopy equivalences of ${\mathbb{Z}}[G]$--modules. The converse holds for  $\pi_1$--injective triples $(\{\kappa_{i}: G_i \rightarrow G\}_{i \in J}, \omega, \mu)$ with $G$ finitely presentable.
\begin{theorem}[{[REALISATION]}]\label{realthm}
There is a ${\rm{PD}}^3$--pair realising the $\pi_1$--injective triple $(\{\kappa_{i}: G_i \rightarrow G\}_{i \in J}, \omega, \mu)$ if and only if $\nu(\mu)$ is a class of homotopy equivalences of ${\mathbb{Z}}[G]$--modules.
\end{theorem}
The case $\partial X = \emptyset$ yields Turaev's Realisation Theorem for ${\rm{PD}}^{3}$--complexes. Precise statements of the results on realisation are contained in Section 4. As the assumption of $\pi_1$--injectivity is indispensable for the method used, the question remains whether the realisation theorem holds without it.

\smallskip

Finally, in Section 5, we use the Classification and Realisation Theorems to show that a ${\rm{PD}}^3$--pair decomposes as connected sum if and only if its $\pi_1$--system decomposes as free product. There are two distinct notions of connected sum for ${\rm{PD}}^3$--pairs reflecting the situation of manifolds with boundary.  We introduce the interior connected sum of two ${\rm{PD}}^3$--pairs as well as the connected sum along boundary components.

If the ${\rm{PD}}^3$--pair $(X, \partial X)$ is orientedly homotopy equivalent to the interior connected sum of two ${\rm{PD}}^3$--pairs $(X_k, \partial X_k)$ with $\pi_1$--systems $\{ \kappa_{k \ell}: G_{k \ell} \rightarrow G_k \}_{\ell \in J_k}$ for $k = 1, 2$, then the $\pi_1$--system of $(X, \partial X)$ is isomorphic to $\{ \iota_k \circ \kappa_{k \ell}: G_{k \ell} \rightarrow G_1 \ast G_2 \}_{\ell \in J_k, k = 1, 2}$, where $\iota_k: G_k \rightarrow G_1 \ast G_2$ denotes the inclusion of the factor in the free product of groups for $ k = 1, 2$. We then say that the $\pi_1$--system of $(X, \partial X)$ \emph{decomposes as free product}. The Classification Theorem  and the Realisation Theorem allow us to show that the converse holds for finitely dominated ${\rm{PD}}^3$--pairs with aspherical boundaries in the case of $\pi_1$--injectivity. 
\begin{theorem}[{[DECOMPOSITION I]}]\label{decomp1}
Let $(X, \partial X)$ be a finitely dominated ${\rm{PD}}^3$--pair with aspherical boundary components. Then $(X, \partial X)$ decomposes as interior connected sum of two $\pi_1$--injective ${\rm{PD}}^3$--pairs if and only if its $\pi_1$--system decomposes as free product of two injective $\pi_1$--systems. 
\end{theorem}

If the ${\rm{PD}}^3$--pair $(X, \partial X)$ is orientedly homotopy equivalent to the boundary connected sum of two ${\rm{PD}}^3$--pairs $(X_k, \partial X_k)$ with $\pi_1$--systems $\{ \kappa_{k \ell}: G_{k \ell} \rightarrow G_k \}_{\ell \in I_k}, k = 1, 2$, along the boundary components $\partial X_{1 \ell_1}$ and $\partial X_{2 \ell_2}$, the $\pi_1$--system of $(X, \partial X)$ is isomorphic to $\{ \kappa: K \rightarrow G_1 \ast G_2, \ \iota_k \circ \kappa_{k \ell}: G_{k \ell} \rightarrow G_1 \ast G_2 \}_{ \ell \in I_k, \ell \neq \ell_k, k = 1, 2},$ where $K := \pi_1(\partial X_{1,\ell_1} \sharp \partial X_{2,\ell_2}; \ast)$ and $\kappa: K \rightarrow G_1 \ast G_2$ is induced by the inclusion of the connected sum of the boundary components $\partial X_{1 \ell_1}$ and $\partial X_{2 \ell_2}$ in the connected sum of the pair. We then say that the $\pi_1$--system of $(X, \partial X)$ \emph{decomposes as free product along $G_{1 \ell_1}$ and $G_{2 \ell_2}$}. As for the interior connected sum, the converse holds for finitely dominated 
${\rm{PD}}^3$--pairs with non--empty aspherical boundaries in the case of $\pi_1$--injectivity.
\begin{theorem}[{[DECOMPOSITION II]}]\label{decomp2}
Let $(X, \partial X)$ be a finitely dominated ${\rm{PD}}^3$--pair  with non--empty aspherical boundary components. Then $(X, \partial X)$ decomposes as boundary connected sum of two $\pi_1$--injective ${\rm{PD}}^3$--pairs $(X_k, \partial X_k), k = 1, 2,$ along $\partial X_{1 \ell_1}$ and $\partial X_{2 \ell_2}$, if and only if its $\pi_1$--system decomposes as free product of two injective $\pi_1$--systems along 
$\pi_1(\partial X_{1 \ell_1}; \ast)$ and $\pi_1(\partial X_{2 \ell_2}; \ast).$
\end{theorem}

A recent paper by Wall \cite{Wall5} contains a different formulation of the first decomposition theorem.

Given an orientable ${\rm{PD}}^3$--pair $(X, \partial X)$, where $\partial X$ is an aspherical $2$--manifold, Crisp's algebraic loop theorem  \cite{Crisp2} allowed him to construct a $\pi_1$--injective ${\rm{PD}}^3$--pair $(\hat X, \partial \hat X)$ with the same fundamental group. Using the decomposition theorem for interior connected sums,  Crisp showed that $(\hat X, \partial \hat X)$ is homotopy equivalent to the connected sum of a finite number of aspherical ${\rm{PD}}^3$--pairs and a  ${\rm{PD}}^3$--complex with virtually free fundamental group. The result holds for the non--orientable case with the addition that the connected summand with virtually free fundamental group may be a ${\rm{PD}}^3$--pair whose boundary is a disjoint union of copies of ${\mathbb{R}} P^2$.

\section{The Relative Twisted Cap Products and Algebraic Sums}
Let ${{\Lambda}}:= {\mathbb{Z}}[G]$ be the integral group ring of the group $G$ and ${\rm{aug}}: {{\Lambda}} \rightarrow {\mathbb{Z}}$ the
\emph{augmentation homomorphism}, defined by ${\rm{aug}}(g) := 1$ for all $g \in G$. Its kernel is the \emph{augmentation ideal}, $I$. A cohomology class, $\omega \in {\rm{Hom}}(G; {\mathbb{Z}}/2{\mathbb{Z}})$, determines a homomorphism from $G$ to ${\mathbb{Z}} / 2 {\mathbb{Z}} = \{0,1\}$, which, in turn, gives rise to the anti--isomorphism $\overline{\phantom{x}}: {{\Lambda}} \rightarrow {{\Lambda}}$ defined by $\overline g := (-1)^{\omega(g)} g^{-1}$ for every $g \in G$. This anti--isomorphism allows us to associate a left ${{\Lambda}}$--module with every right ${{\Lambda}}$--module and vice versa. Namely, given a right ${{\Lambda}}$--module $A$  we define a left action on the set underlying $A$ by ${{\lambda}}.a := a.\overline{{\lambda}}$ for every $a \in A$ and ${{\lambda}} \in {{\Lambda}}$. Proceeding analogously for a left ${{\Lambda}}$--module $B$, we obtain a left ${{\Lambda}}$--module $^{\omega}\!\!A$ and a right ${{\Lambda}}$--module $B^{\omega}$.

Take chain complexes $C$ and $D$ of left ${{\Lambda}}$--modules and a chain complex $E$ of right ${{\Lambda}}$--modules. Define $E \otimes_{{{\Lambda}}} D$ and ${\rm{Hom}}_{\Lambda}(C,D)$ as usual with the conventions
\begin{gather*}
\partial^{E \otimes_{{{\Lambda}}} D} |_{E_{i} \otimes_{{{\Lambda}}} D_{j}} = \partial^{E}
\otimes {\mbox{id}} + (-1)^{i} {\mbox{id}} \otimes \partial^{D}; \\
\partial \big( (f_{i})_{i \in {\mathbb{Z}}} \big)|_{{\rm{Hom}}_{\Lambda}(C,D)_n}:= (\partial^{D} f_{i})_{i \in
{\mathbb{Z}}} - \big( (-1)^{n} f_{i} \partial^{C} \big)_{i \in {\mathbb{Z}}}.
\end{gather*}
Treating right ${{\Lambda}}$--modules $A$ and left ${{\Lambda}}$--modules $B$ as chain complexes concentrated at
level zero, we define
\[ {\mbox{H}}_{n}(C;A) := {\mbox{H}}_{n}(A {\otimes_{\Lambda}} C); \quad
{\mbox{H}}^{k}(C;B) := {\mbox{H}}_{-k}({\rm{Hom}}_{\Lambda}(C,B)).\]
The group $G$ acts diagonally on $C {\otimes}D$ via $g(c_i \otimes d_j) = g c_i \otimes g d_j$ for $g \in G$ and $c_i \otimes d_j \in C_i \otimes D_j$. On the level of chains, the twisted slant
operation is the chain map
\begin{gather*}
/ : \big({\rm{Hom}}_{\Lambda}(C, B)\big)_{-k} {\otimes}\big( {\mathbb{Z}}^{\omega} {\otimes_{\Lambda}} (C {\otimes}D)
\big)_{n}
\longrightarrow B^{\omega} {\otimes_{\Lambda}} D_{n-k}, \\
\varphi / (z {\otimes}\sum_{i+j=n} c_{i} {\otimes}d_{j}) := z \varphi(c_{k}) \otimes
d_{n-k}.
\end{gather*}
Passing to homology and composing with the homomorphism 
\[ \alpha: {\mbox{H}} C {\otimes}{\mbox{H}} D \rightarrow {\mbox{H}} (C {\otimes}D), [x] {\otimes}[y] \mapsto [x {\otimes}y],\]
yields the twisted slant operation
\[
/ : {\mbox{H}}^{k}(C; B) {\otimes}{\mbox{H}}_{n}(C {\otimes}D; {\mathbb{Z}}^{\omega}) \rightarrow
{\mbox{H}}_{n-k}(D; B^{\omega}).
\]

Let $\varepsilon: P \twoheadrightarrow {\mathbb{Z}}$ be an augmented chain complex of left ${{\Lambda}}$--modules with equivariant diagonal $\Delta: P \rightarrow P {\otimes}P$, that is, $(\varepsilon {\otimes}{\rm{id}}) \Delta(c) = ({\rm{id}} {\otimes}\varepsilon) \Delta(c) = c$, and let $Q$ be a subcomplex of $P$ such that the inclusion $\iota: Q \rightarrowtail P$ is a map of augmented chain complexes with equivariant diagonal, that is, $\Delta \circ \iota = (\iota {\otimes}\iota) \circ \Delta$. Then $(P, Q)$ is called a \emph{geometric chain pair}, 
$Q \stackrel{\iota}{\rightarrowtail} P \stackrel{\pi}{\twoheadrightarrow } D:= P/Q$ a \emph{short exact sequence of augmented chain complexes with compatible diagonals} and the chain map
\[\Delta_{{\rm{rel}}}: D \longrightarrow P {\otimes}D, d \mapsto (\rm{id} {\otimes}\pi)(\Delta(p)), \quad {\rm{where}} \quad d = \pi(p), \]
is called the \emph{relative equivariant diagonal}. Composing the relative diagonal and the 
twisted slant operation, we obtain the \emph{relative twisted cap product}
\[ \cap: {\rm{Hom}}_{\Lambda}(P, B)_{-k} {\otimes}({\mathbb{Z}}^{\omega} {\otimes_{\Lambda}} D)_{n} \rightarrow (B^{\omega}
{\otimes_{\Lambda}} D)_{n-k}, (\varphi, (z {\otimes}d)) \mapsto \varphi / (z {\otimes}\Delta_{{\rm{rel}}}(d))\]
for every left ${{\Lambda}}$--module $B$. Passing to homology and composing with $\alpha$, we obtain the relative twisted cap product
\[ \cap: {\mbox{H}}^{k}(P; B) {\otimes}{\mbox{H}}_{n}(D;{\mathbb{Z}}^{\omega}) \rightarrow {\mbox{H}}_{n-k}(D; B^{\omega}). \]
Similarly, the relative diagonal $\Delta'_{\rm{rel}}: D \longrightarrow D {\otimes}P, d \mapsto (\pi {\otimes}\rm{id})(\Delta(p))$,  where $d = \pi(p)$, yields the relative twisted cap product
\[\cap: {\mbox{H}}^{k}(D; B) {\otimes}{\mbox{H}}_{n}(D;{\mathbb{Z}}^{\omega}) \rightarrow {\mbox{H}}_{n-k}(P; B^{\omega}).\]
Note that the relative cap products reduce to the absolute cap product when $Q$ is trivial.

If  $M$ is a ${{\Lambda}}$--bimodule, then ${{\Lambda}}$ acts on the left on $M {\otimes_{\Lambda}} B$ via 
$\lambda.(m {\otimes}b) := (\lambda.m) {\otimes}b$ and on the right on
${\rm{Hom}}_{\Lambda}(B, M)$ via $(\varphi.\lambda)(b) := \varphi(b).\lambda$ for all left 
${{\Lambda}}$--modules $B, \lambda \in {{\Lambda}}, b \in B, m \in M$ and 
$\varphi \in {\rm{Hom}}_{\Lambda}(B, M)$. In particular, $B^{\ast} := ^{\omega}\!{{\rm{Hom}}_{\Lambda}(B, {{\Lambda}})}$ 
is a left ${{\Lambda}}$--module. Any left ${{\Lambda}}$--module $A$ gives rise to the functors 
$A^{\omega} {\otimes_{\Lambda}} -$ and ${\rm{Hom}}_{\Lambda}({^{\omega}{\rm{Hom}}_{\Lambda}( - , {{{\Lambda}}})}, A))$ from the category
$_{{{\Lambda}}}\mathcal M$ of left ${{\Lambda}}$--modules to the
category $\mathcal Ab$ of abelian groups and there is a natural 
transformation 
\begin{eqnarray}\label{eta}
\eta_{B}: A^{\omega} {\otimes_{\Lambda}} B \longrightarrow
{\rm{Hom}}_{\Lambda}(B^{\ast}, A)
\end{eqnarray}
given by
\[ \eta_{B}(a {\otimes}b): \ ^{\omega}\!{\rm{Hom}}_{\Lambda}(B, {{{\Lambda}}}) \longrightarrow A, \quad \varphi 
\longmapsto \overline{\varphi(b)}a \]
for every left ${{\Lambda}}$--module $B$. When we restrict the two functors 
to the category of finitely generated free left ${{\Lambda}}$--modules, the natural transformation 
$\eta$ becomes a natural equivalence. For ${{\Lambda}}$--bimodules $M$ we may view
$^{\omega}\!{{\rm{Hom}}_{\Lambda}( ^{\omega}\!{\rm{Hom}}_{\Lambda}( -, {{\Lambda}}) , M)}$ 
and $^{\omega}\!M^{\omega} {\otimes_{\Lambda}} - $ as functors from the
category of left ${{\Lambda}}$--modules to itself and in this case the natural transformation $\eta$ 
respects the additional left ${{\Lambda}}$--module structure. Identifying the left ${{\Lambda}}$--module $B$ with $^{\omega}\!{{\Lambda}}^{\omega} {\otimes_{\Lambda}} B$ for $M = {{\Lambda}}$, the natural equivalence
$\eta$ becomes the evaluation homomorphism from $B$ to its double dual
\[(B^{\ast})^{\ast} = ^{\omega}\!{{\rm{Hom}}_{\Lambda}( ^{\omega}\!{\rm{Hom}}_{\Lambda}( -, {{\Lambda}}) , {{\Lambda}})}.\]

Let $Q \rightarrowtail P \twoheadrightarrow D$ be a short exact
sequence of augmented chain complexes with compatible diagonals,
so that $(P, Q)$ is a geometric chain pair. Then the chain map given by taking 
the cap product with a cycle $1 {\otimes}x \in {\mathbb{Z}}^{\omega} {\otimes_{\Lambda}} D_{n}$ is almost chain 
homotopic to its dual. 
\begin{lemma}\label{capastcom}
Let $1 {\otimes}x \in {\mathbb{Z}}^{\omega} {\otimes_{\Lambda}} D_{n}$ be a cycle. Then the diagram
\[ \xymatrix{
D_{k}^{\ast} \ar[rrr]^{\theta} \ar[d]_{\cap 1 {\otimes}x} &&&
( ^{\omega}\!{{\Lambda}}^{\omega}\! {\otimes_{\Lambda}} D_{k})^{\ast} \ar[d]^{(\cap 1 {\otimes}x)^{\ast}} \\
^{\omega}\!{{\Lambda}}^{\omega}\! {\otimes_{\Lambda}} P_{n-k}
\ar[rrr]_{\eta_{P_{n-k}}} &&& (P_{n-k}^{\ast})^{\ast}} \]
commutes up to chain homotopy, where the isomorphism $\theta$
is given by $\theta(\varphi)(\lambda {\otimes}d) := \overline{\lambda}\varphi(d)$
for every $\varphi \in D_k^{\ast}, \lambda \in {{\Lambda}}$ and $d \in D_k$.
\end{lemma}
\begin{proof}
Use that any two diagonals are chain homotopic \cite{Bleile}.
\end{proof}
The short exact sequence $Q \stackrel{\iota}{\rightarrowtail} P \stackrel{\pi}
{\twoheadrightarrow} D$ of augmented chain
complexes of free ${{\Lambda}}$--modules with compatible diagonals
splits and stays split short exact when we tensor or apply the
${\rm{Hom}}_{\Lambda}$--functor. Given a right ${{\Lambda}}$--module $M$, denote the connecting
homomorphisms of 
\begin{gather*}
{\mathbb{Z}}^{\omega} {\otimes_{\Lambda}} Q \rightarrowtail {\mathbb{Z}}^{\omega} {\otimes_{\Lambda}} P
\twoheadrightarrow {\mathbb{Z}}^{\omega} {\otimes_{\Lambda}} D, \\
M {\otimes_{\Lambda}} Q \rightarrowtail M {\otimes_{\Lambda}} P \twoheadrightarrow M {\otimes_{\Lambda}} D \quad {\rm{and}} \quad \\
^{\omega}\!{\rm{Hom}}_{\Lambda}(D, ^{\omega}\!\!{M}) \rightarrowtail ^{\omega}\!\!{\rm{Hom}}_{\Lambda}(P, ^{\omega}\!\!{M}) \twoheadrightarrow^{\omega}\!\!{\rm{Hom}}_{\Lambda}(Q, ^{\omega}\!\!{M})
\end{gather*}
by $\delta_{\ast}$, $\delta_{\ast}'$ and $\delta^{\ast}$ respectively. Standard arguments \cite{Bleile}
show

\begin{theorem}[Cap Product Ladder]\label{ladder}
For all $y \in {\mbox{H}}_{n}(D;{\mathbb{Z}}^{\omega})$, the diagram
\[ \xymatrix{
\cdots \ar[r] & {\mbox{H}}^{r}(D; ^{\omega}\!\!M) \ar[r]^{\pi^{\ast}} \ar[d]^{\cap y} &
{\mbox{H}}^{r}(P; ^{\omega}\!\!M) \ar[r]^{\iota^{\ast}} \ar[d]^{\cap y} &
{\mbox{H}}^{r}(Q; ^{\omega}\!\!M) \ar[d]^{\cap \delta_{\ast} y} \\
\cdots \ar[r] & {\mbox{H}}_{n-r}(P;M) \ar[r] & {\mbox{H}}_{n-r}(D;M) \ar[r]^-{\delta_{\ast}'} & {\mbox{H}}_{n-r-1}(Q;M) \\
\ar[r]^-{\delta^{\ast}} & {\mbox{H}}^{r+1}(D; ^{\omega}\!\!M) \ar[r] \ar[d]^{\cap y} & \cdots \\
\ar[r] & {\mbox{H}}_{n-r-1}(P;M) \ar[r] & \cdots} \]
commutes, up to sign.
\end{theorem}
The geometric chain pair $(P, Q)$ is called a \emph{Poincar\'{e} chain pair} if there is an element 
$\nu \in H_n(P; Q; {\mathbb{Z}}^{\omega})$ of infinite order such that 
\[ \cap \nu: {\mbox{H}}^k(P; ^{\omega}\!\!M) {\otimes}{\mbox{H}}_{n}(D;{\mathbb{Z}}^{\omega}) \rightarrow {\mbox{H}}_{n-k}(D; M)\]
is an isomorphism of abelian groups for every right ${{\Lambda}}$--module $M$ and every $k \in {\mathbb{Z}}$. We then call $\nu$ the \emph{fundamental class} of $(P, Q)$. For the definition of connected sums of $\rm{PD}^n$--pairs we need a generalisation of Browder's result \cite{Browder} concerning cap products and sums of Poincar\'e chain pairs to the case of non--trivial fundamental groups.

Suppose $B = B_1 + B_2, B_0 = B_1 \cap B_2, A \subseteq B,  A_i = B_i \cap A$ for $i = 1, 2,$ and $A_0 = A_1 \cap A_2$ are chain complexes of free ${{\Lambda}}$--modules such that all pairs arising are geometric chain pairs. Let $\partial_0$ be the connecting homomorphism of the short exact sequence
\[
\xymatrix{
{\mathbb{Z}}^{\omega} {\otimes_{\Lambda}}(B_0/A_0) \ \ar@{>->}[r] & {\mathbb{Z}}^{\omega} {\otimes_{\Lambda}} (B_1/A_1 \oplus B_2/A_2)
\ar@{->>}[r] & {\mathbb{Z}}^{\omega} {\otimes_{\Lambda}} (B/A),}\]
and, for $i = 1, 2,$ let $\eta_i$ denote the ${{\Lambda}}$--morphism rendering the diagram
\[
\xymatrix{
H_q(B, A; {\mathbb{Z}}^{\omega}) \ar[rr]^{\eta_i} \ar[drr] && H_q(B_1, B_0 + A_1; {\mathbb{Z}}^{\omega}) \ar@{>->>}[d] \\
&& H_q(B, B_2 + A; {\mathbb{Z}}^{\omega})}\]
commutative. Standard arguments \cite{Bleile} show
\begin{theorem}\label{sumthm}
Any two of the following conditions imply the third:
\begin{enumerate}
\item $(B, A)$ is a Poincar\'{e} chain pair with fundamental class $\nu \in H_n(B, A; {\mathbb{Z}}^{\omega})$;
\item $(B_0, A_0)$ is a Poincar\'{e} chain pair with fundamental class $\partial_0 \nu$ an element of $H_{n-1}(B_0, A_0; {\mathbb{Z}}^{\omega})$;
\item for $i = 1, 2, \ (B_i, B_0 + A_i)$ is a Poincar\'{e} chain pair with fundamental class $\eta_i \nu \in H_n(B_i, B_0 + A_i; {\mathbb{Z}}^{\omega})$.
\end{enumerate}
\end{theorem}
\section{Homotopy Classification of $\rm{PD}^3$--Pairs}
Given a CW--complex  $X$, fix a base point $\ast$ in $X$, put $G :=
\pi_{1}(X; \ast)$ and ${{\Lambda}} := {\mathbb{Z}}[G]$. Let $p: \tilde X \rightarrow X$ denote the universal 
covering of $X$ and $C(X)$ the cellular chain complex of $\tilde X$, viewed as a 
complex of left ${{\Lambda}}$--modules. Then, for a right ${{\Lambda}}$--module $A$ and a left 
${{\Lambda}}$--module $B$,
\[ {\mbox{H}}_{r}(X;A) := {\mbox{H}}_{r}(C(X);A) \quad {\rm{and}} \quad
{\mbox{H}}^{r}(X;B) := {\mbox{H}}^{r}(C(X);B). \]
A \emph{connected ${\rm{PD}}^{n}$--complex} is a triple $(X, \omega, [X])$, 
where $X$ is a connected CW--complex, $\omega \in {\mbox{H}}^{1}(X;{\mathbb{Z}} / 2 {\mathbb{Z}})$ and 
$[X] \in {\mbox{H}}_{n}(X;{\mathbb{Z}}^{\omega})$ such that
\begin{equation}\label{capwithfund}
    \cap [X]: {\mbox{H}}^{r}(X;B) \rightarrow {\mbox{H}}_{n-r}(X; B^{\omega});
    \quad \alpha \mapsto \alpha \cap [X]
\end{equation}
is an isomorphism of abelian groups for every $r \in {\mathbb{Z}}$ and every
left ${{\Lambda}}$--module $B$. We often write $X$ for the triple $(X, \omega, [X])$
and call $\omega$ the \emph{orientation character} and $[X]$ the \emph{fundamental
class} of $X$. A ${\rm{PD}}^{n}$--complex $(X, \omega, [X])$ is a finite disjoint union
of connected ${\rm{PD}}^{n}$--complexes $(X_{i}, \omega_{i}, [X]_{i})$,
$i \in J$, where ${{\Lambda}} = \bigoplus_{i \in J} {\mathbb{Z}} [G_{i}], \ \omega = (\omega_{i})_{i \in J}$ 
and $[X] = ([X_{i}])_{i \in J}$. Note that it is enough to demand
that~(\ref{capwithfund}) is an isomorphism for $B = {{\Lambda}}$ \cite{Wall1}.

Every $n$--dimensional manifold is homotopy equivalent to a $CW$--complex and 
thus determines a ${\rm{PD}}^{n}$--complex. But not every ${\rm{PD}}^3$--complex 
is homotopy equivalent to a $3$--manifold. Wall \cite{Wall1} showed 
that the class of finite groups which are fundamental groups of 
${\rm{PD}}^3$--complexes coincides with the class of finite groups with periodic 
cohomology of period $4$. By Milnor's results \cite{Milnor1}, some of these groups 
are not realisable as fundamental groups of $3$--manifolds, the simplest such group 
being the permutation group $S_3$. Swan \cite{Swan} explicitly constructed
a ${\rm{PD}}^3$--complex $K$ with fundamental group $S_3$.

Given a pair, $(X, \partial X)$, of CW--complexes, let $C(\partial X)$ denote the 
subcomplex of $C(X)$ generated by the cells lying above $\partial X$ and put
$C(X,\partial X) := C(X) /C(\partial X)$. Then $C(X, \partial X)$ is called the 
\emph{relative cellular complex} and 
\[ C(\partial X) \rightarrowtail C(X) \twoheadrightarrow C(X,\partial X)\]
the \emph{short exact sequence of cellular chain complexes} of the pair $(X,\partial X)$. 
We define
\[ {\mbox{H}}_{r}(X,\partial X;A) := {\mbox{H}}_{r}(C(X,\partial X);A) \quad {\rm{and}} \quad
{\mbox{H}}^{r}(X,\partial X;B) := {\mbox{H}}^{r}(C(X,\partial X);B)\]
for every right ${{\Lambda}}$--module $A$ and every left ${{\Lambda}}$--module $B$, and denote the 
connecting homomorphism of ${\mathbb{Z}}^{\omega} {\otimes_{\Lambda}} C(\partial X) \rightarrowtail 
{\mathbb{Z}}^{\omega} {\otimes_{\Lambda}} C(X) \twoheadrightarrow {\mathbb{Z}}^{\omega} {\otimes_{\Lambda}} C(X,\partial X)$ by 
$\delta_{\ast}$.
\begin{definition}
The quadruple $(X, \partial X, \omega_{X}, [X,\partial X])$ is a \emph{connected 
${\rm{PD}}^{n}$--pair} if $(X, \partial X)$ is a pair of CW--complexes with $X$ connected,
$(\partial X, \omega_{\partial X}, [\partial X])$ is a ${\rm{PD}}^{n-1}$--complex, 
where $\omega_{X}$ induces $\omega_{\partial X_{i}}$ on the connected
components $\partial X_{i}$ of $\partial X$, and $[X,\partial X]$ is an element of 
${\rm{H}}_{n}(X, \partial X; {\mathbb{Z}}^{\omega})$ with $\delta_{\ast}[X,\partial X] = [\partial X]$,
such that
\begin{equation}\label{capwithpairfund}
\cap [X,\partial X]: {\rm{H}}^{r}(X; B) \rightarrow {\rm{H}}_{n-r}(X, \partial X; B^{\omega}), 
\alpha \mapsto \alpha \cap [X,\partial X]
\end{equation}
is an isomorphism for every left ${{\Lambda}}$--module $B$ and every $r \in {\mathbb{Z}}$.
We often write $(X, \partial X)$ for the ${\rm{PD}}^{n}$--pair and call
$\omega_{X}$ the \emph{orientation character} and $[X,\partial X]$ the
\emph{fundamental class} of $(X, \partial X)$.
\end{definition}
Again it is enough to demand that~(\ref{capwithpairfund}) is an
isomorphism for $B = {{\Lambda}}$ \cite{Wall1}.
As for ${\rm{PD}}^{n}$--complexes, we may define a general
${\rm{PD}}^{n}$--pair as a finite disjoint union of connected
${\rm{PD}}^{n}$--pairs. Every manifold with boundary determines a
${\rm{PD}}^{n}$--pair. We obtain examples of ${\rm{PD}}^{n}$--pairs 
which are not homotopy equivalent to a manifold with boundary by
taking a ${\rm{PD}}^{n}$--complex which is not homotopy equivalent
to a manifold and removing the interior of an $n$--cell.

Poincar\'e duality manifests itself on the level of chains, namely, for a ${\rm{PD}}^{n}$--pair, 
$(X, \partial X)$, the cap product with a representative of the fundamental 
class defines a chain map of degree $n$ which is a chain homotopy equivalence.
\begin{definition}
The ${\rm{PD}}^{n}$--pairs $(X_{1}, \partial X_{1})$ and
$(X_{2}, \partial X_{2})$ are {\emph{orientedly homotopy equivalent}}
if there is a homotopy equivalence $f: (X_{1}, \partial X_{1}) \rightarrow
(X_{2}, \partial X_{2})$ of pairs with
$f^{\ast}(\omega_{X_{2}}) = \omega_{X_{1}}$ and $f_{\ast}([X_{1},\partial X_{1
}])
= [X_{2},\partial X_{2}]$.
\end{definition}
A family, $\{\kappa_{i}: G_{i} \rightarrow G\}_{i \in J}$, of group homomorphisms is 
also called a \emph{$\pi_1$--system} and an \emph{Eilenberg--{Mac\,Lane} pair} of
type $K(\{\kappa_{i}: G_{i} \rightarrow G\}_{i \in J};1)$ is a pair $(X, \partial X)$
such that $X$ is an Eilenberg--{Mac\,Lane} complex of type $K(G;1)$, the
connected components $\{\partial X_{i}\}_{i \in J}$ of $\partial X$ are
Eilenberg--{Mac\,Lane} complexes of type $K(G_{i};1)$ and there is an isomorphism
\[ (\varphi, \{\varphi_{i}\}_{i \in J}): \{ \rho_{i}: G_{i} \rightarrow G\}_{i \in J} \rightarrow 
\{ \kappa_{i}: \pi_1(\partial X_{i}, \ast) \rightarrow \pi_1(X, \ast) \}_{i \in J}\]
of $\pi_1$--systems, that is, for each $i \in J$ the diagram
\[
\xymatrix{
G_{i} \ar[r]^-{\varphi_{i}} \ar[d]_{\rho_{i}}
& \pi(\partial X_{i}, \ast) \ar[d]_{\kappa_{i}} \\
G \ar[r]^-{\varphi} & \pi_1(X, \ast)}
\]
commutes up to conjugacy. Note that we do not require the homomorphisms
$\kappa_{i}$ to be injective as in the standard definition given by Bieri--Eckmann 
~\cite{Bieri-Eckmann}. For any $\pi_1$--system, $\{\kappa_{i}: G_{i} \rightarrow G\}_{i \in J}$, 
there is an Eilenberg--{Mac\,Lane} pair of type $K(\{ \kappa_i: G_{i} \rightarrow G \}_{i \in J};1)$, 
which is uniquely determined up to homotopy equivalence of pairs. Further, for any pair, 
$(X, \partial X)$, of CW--complexes there is map of pairs
\[ c_{X}: (X, \partial X) \longrightarrow K(G,\{G_{i}\}_{i \in J};1), \]
called the \emph{classifying map}, which is uniquely determined up to homotopy 
of pairs and induces an isomorphism of $\pi_1$--systems.

An oriented homotopy equivalence of ${\rm{PD}}^{n}$--pairs, $(X,\partial X)$ 
and $(Y,\partial Y)$, induces an isomorphism 
\[ (\varphi, \{\varphi_{i}\}_{i \in J}):
\{ \kappa_{i}: \pi_1(\partial X_{i}, \ast) \rightarrow \pi_1(X, \ast) \}_{i \in J} \rightarrow 
\{ \rho_{i}: \pi_1(\partial Y_{i}, \ast) \rightarrow \pi_1(Y, \ast) \}_{i \in J}\]
of their $\pi_{1}$--systems such that
\begin{equation}\label{grpsysisoplus}
\varphi^{\ast}(\omega_{Y}) = \omega_{X} \quad {\rm{and}} \quad
\varphi_{\ast}(c_{X \ast}([X, \partial X])) = c_{Y \ast}([Y, \partial Y]),
\end{equation}
where $c_{X}$ and $c_{Y}$ are classifying maps  and $\{ \partial X_i \}$ and 
$\{ \partial Y_i \}$ are the connected components of $\partial X$ and $\partial Y$ 
respectively. The Classification Theorem states that there is a converse for $n=3$.

Now put $\mu_X := c_{X \ast}([X, \partial X])$ and restrict attention to the case $n=3$. 

As the homology and cohomology sequences of any Eilenberg--{Mac\,Lane} pair
$K(\{\kappa_{i}: G_{i} \rightarrow G\}_{i \in J};1)$ are isomorphic to those of the group
pair $(G, \{G_i\}_{i \in J})$ \cite{Bieri-Eckmann}, we may identify $c_{X \ast}([X, \partial X]$ 
and $\omega_X$ with their images in  ${\mbox{H}}^{1}(G;{\mathbb{Z}} / 2 {\mathbb{Z}})$ and $\mbox{H}_3(G, \{G_i\}_{i \in J}; {\mathbb{Z}}^{\omega})$, respectively.
\begin{definition}
The triple $(\{ \kappa_{i} \}_{i \in J}, \omega_X, \mu_X)$ is the \emph{fundamental triple}
of the ${\rm{PD}}^{3}$--pair $(X,\partial X)$.

Two triples $(\{ \kappa_{i} \}_{i \in J}, \omega, \mu)$ and 
$(\{ \kappa'_{i} \}_{i \in J}, \omega', \mu')$ are \emph{isomorphic} if there is an 
isomorphism $(\varphi, \{\varphi_{i}\}_{i \in J}): \{ \kappa_{i} \}_{i \in J} \rightarrow 
\{ \kappa'_{i} \}_{i \in J}$ of $\pi_1$--systems such that $\varphi^{\ast}(\omega') = \omega$ and $\varphi_{\ast}(\mu) = \mu'$.
\end{definition}
Note that the fundamental triple of a ${\rm{PD}}^{3}$--pair is uniquely determined up to 
isomorphism of triples.
\begin{proof}[Proof of Classification]
The proof generalises Turaev's alternative proof of Hendriks' result for the absolute
case.

Given two ${\rm{PD}}^{3}$--pairs, $(X,\partial X)$ and $(Y,\partial Y)$, with 
aspherical boundaries and isomorphic fundamental triples, we use homological algebra and 
obstruction theory to construct a map $f: (X,\partial X) \rightarrow (Y, \partial Y)$ of pairs such 
that
\begin{itemize}
\item[(i)] $f$ induces a homotopy equivalence $\hat f: \partial X
\rightarrow \partial Y$; 
\item[(ii)] $f_{\ast}: \pi_{1}(X,\ast) \rightarrow \pi_{1}(Y,\ast)$ is
an isomorphism respecting the orientation character; 
\item[(iii)] $f$ has degree one, that is, $f_{\ast}([X,\partial X]) =
[Y,\partial Y]$.
\end{itemize}
Then $f$ is a homotopy equivalence of ${\rm{PD}}^{3}$--pairs 
by Poincar{\'e} duality, the Five Lemma and Whitehead's Theorem.

It is not difficult to see that $X$ is homologically two--dimensional for any 
${\rm{PD}}^{3}$--pair, $(X, \partial X)$. By definition, the boundary $\partial X$ is
a ${\rm{PD}}^{2}$--complex and Eckmann, M\"{u}ller and Linnell 
(\cite{Eckmann-Mueller} and \cite{Eckmann-Linnell}) showed that every
connected ${\rm{PD}}^{2}$--complex is homotopy equivalent to a closed
surface. If the boundary is not empty, we may apply
the mapping cylinder construction and hence assume that
the components $\{ \partial X_{i} \}_{i \in J}$ of $\partial X$ are closed surfaces 
which are collared in $X$, the collar being a $3$--manifold with boundary. 
Splitting off a $3$--cell from the collar of one of the boundary components, we obtain 
$X = X' \cup_{g} e^{3}$ where $g: S^{2} \rightarrow X'$.

By definition \cite{Wall2}, the connected CW--complex $X$ satisfies $\rm{D}_n$ if ${\mbox{H}}_i(\tilde X) = 0$ for $i >n$ and ${\mbox{H}}^{n+1}(X; B) = 0$ for all left ${{\Lambda}}$--modules $B$. Poincar\'e duality implies that $X$ satisfies $\rm{D}_2$ for every ${\rm{PD}}^3$--pair $(X, \partial X)$. We obtain
\cite{Bleile}
\begin{lemma}\label{splitoffe3}
Let $(X, \partial X)$ be a ${\rm{PD}}^3$--pair. 
Then $X = X' \cup_{g} e^{3}$ where $g: S^{2} \rightarrow X'$
and $X'$ satisfies $\rm{D}_2$ and is thus homologically two--dimensional.
\end{lemma}
Take a connected ${\rm{PD}}^{3}$--pair, $(X,\partial X)$, with aspherical boundary components and
$X = X' \cup_{g} e^{3}$ as above. Without loss of generality, we may assume 
that $\partial X \cap e^{3} = \emptyset$. We denote the subcomplex of $C(X)$
generated by the cells lying in $\tilde X' := p^{-1}(X')$ by $C(X')$. By Theorem [E] in \cite{Wall2}, we may assume that $X$ and hence $X'$ are geometrically $3$--dimensional. As $X'$ is
homologically $2$--dimensional, $C(X')$ is a chain complex of free left ${{\Lambda}}$--modules with 
${\mbox{H}}^k(C(X');B) = 0$ for $k \geq 3$ and for every (left) ${{\Lambda}}$--module $B$. Hence Lemma 3.6 in \cite{Bleile} implies that $C(X')$ decomposes into the direct sum of the two subcomplexes
\[
\begin{array}{cccccccccc}
D: & 0 & \rightarrow & C_{3}(X') & \twoheadrightarrow & {\rm{im}} \partial_{2}' &  & \\ 
&&&&&& \\
E: & 0 & \rightarrow & {{\Lambda}}[e] & \twoheadrightarrow & S & \rightarrow & C_1(X) & \rightarrow & C_{0}(X),
\end{array}
\]
where $D \simeq 0$, $S$ is the cokernel of the boundary operator $\partial_2': C_3(X') \rightarrow C_2(X')$, the chain $e$ corresponds to the $3$--cell attached via $g$ and ${{\Lambda}}[e]$ denotes the 
free ${{\Lambda}}$-module with generator $e$. As before, we may assume that the connected 
components of $\partial X$ are closed surfaces, so that $C_{i}(\partial X) = 0$ for $i > 2$.
\begin{lemma}\label{asphericalcompopensurf}
Let $p: \tilde X \rightarrow X$ be the universal covering. Then the components of
$p^{-1}(\partial X)$ are open surfaces and hence the boundary operator
$C_{2}(\partial X) \rightarrow C_{1}(\partial X)$ is injective.
\end{lemma}
\begin{proof}
Suppose the ${\rm{PD}}^3$--pair $(X, \partial X)$ has $\pi_1$--system 
$\{ \kappa_{i}: G_{i}\rightarrow G \}_{i \in J}$ and connected boundary components 
$\{ \partial X_j \}_{j \in J}$ which are closed aspherical surfaces. Suppose further that $\kappa_ i(G_i)$ is finite in $G$ for some $i \in J$. Passing to $(X^+, \partial X^+)$, where $X^+$ is the orientable covering space of $X$ and $\partial X^+$ is the inverse image of $\partial X$ under the covering map, the image of the fundamental group of each component of $\partial X^+$ over $\partial X_i$ is still finite. Thus we may assume that $(X, \partial X)$ is orientable. 

We replace $X$ by $X' = X \cup_{j \neq i} H_j$, where $(H_j, \partial X_j)$ is a connected ${\rm{PD}}^3$--pair for $j \in J, j \neq i$. By Theorem~\ref{sumthm}, $(X', \partial X_i)$ is a ${\rm{PD}}^3$--pair. Thus we may assume without loss of generality that $\partial X$ is connected.
Poincar{\'e} duality and the Hurewicz homomorphism together with our assumptions yield 
${\mbox{H}}_1(\partial X; {\mathbb{Q}}) = 0$, so that $\partial X = S^2$. Hence, if all components of $\partial X$ are aspherical, all components of $p^{-1}(\partial X)$ are open surfaces.
\end{proof}
Lemma~\ref{asphericalcompopensurf} implies ${\rm{im}} \partial_{2}' \cap C_{2}(\partial X)
 = 0$ and hence the relative chain complex $C(X, \partial X)$ of the pair $(X,\partial
X)$ decomposes into the direct sum of the two subcomplexes
\[
\begin{array}{cccccccccc}
D: & 0 & \rightarrow & C_{3}(X') & \twoheadrightarrow & {\rm{im}} \partial_{2}' & \rm{and} & &\\
&&&&&&&& \\
E: & 0 & \rightarrow & {{\Lambda}}[e] & \rightarrow & \tilde S & \rightarrow & C_{1}(X, \partial X) & \rightarrow & C_{0}(X, \partial X),
\end{array}
\]
where $\tilde S = S / C_{2}(\partial X)$ and $D \simeq 0$.

\smallskip

Let $\{ e^{2}_{m}\}_{m \in M}$ be a collection of two--cells of
$\tilde X$ such that above every two--cell in $X \setminus \partial X$
there lies exactly one cell from this collection. Then the collection
$\{ e_{m}\}_{m \in M}$ of chains represented by these
cells comprises a basis of the ${{\Lambda}}$--module $C_{2}(X, \partial X)$. Suppose
$\partial_{2}(e) = \sum a_{m} e_{m}$, where $\partial_{2}: C_{3}(X, \partial X)
\rightarrow C_{2}(X, \partial X)$ denotes the boundary operator of the relative 
complex and $a_{m} \in {{\Lambda}}$ for $m \in M$.
\begin{lemma}\label{splitoffandmore}
The chain $1 {\otimes}e$ is a relative cycle representing the homology class 
$[X,\partial X]$. Further, $I = {\rm{im}}({\partial}^{E}_{2})^{\ast}$ and is generated by 
$\{ \overline{a_{m}} \}_{m \in M}$.
\end{lemma}

For a proof we refer the reader to \cite{Bleile}.

Given connected ${\rm{PD}}^{3}$--pairs, $(X,\partial X)$ and $(Y,\partial Y)$, with 
aspherical boundary components and isomorphic fundamental triples, we must 
construct a map $f : (X, \partial X) \longrightarrow (Y, \partial Y)$ satisfying (i) -- (iii).
Suppose the isomorphism of fundamental triples is given by
\[ (\varphi, \{\varphi_{i}\}_{i \in J}): \{\kappa_{i}: G_{i} \rightarrow G\}_{i \in J} 
\rightarrow \{ \rho_{i}: H_{i} \rightarrow H\}_{i \in J}.\] 
We may assume without loss of generality that $Y$ is contained in $(K, \partial Y) :=
K( \{ \rho_i: H_{i} \rightarrow H \}_{i \in J};1)$. Let $q: \tilde K \rightarrow K$
be the universal covering, put ${{\Lambda}}_{X} := {\mathbb{Z}}[G]$, ${{\Lambda}}_{Y} := {\mathbb{Z}}[H]$ and let $I_X$ 
and $I_Y$ denote the kernel of ${\rm{aug}}_{X}$ and ${\rm{aug}}_{Y}$ 
respectively. By Lemma~\ref{splitoffandmore}, $X = X' \cup e^{3}$ and $Y = Y' \cup {e'}^{3}$, 
where $X'$ and $Y'$ are homologically two--dimensional, and if $e$ and $e'$ 
are chains representing the $3$--cells $e^{3}$ and ${e'}^{3}$, then $1 {\otimes}e$ 
and $1 {\otimes}e'$ represent $[X, \partial X]$ and $[Y,\partial Y]$ respectively. 

Further, we may assume that the connected components of $\partial X$ and $\partial Y$ 
are closed surfaces, and that there is a homeomorphism 
$\tilde g: \partial X \longrightarrow \partial Y$ which induces the isomorphisms $\varphi_{i}: G_{i} \rightarrow H_{i}$ on the fundamental groups of the connected components of $\partial X$
and $\partial Y$ respectively. It is not difficult to construct a cellular map of pairs
$g': (X', \partial X) \longrightarrow (Y',\partial Y)$ with $g' |_{\partial X} = \tilde g$. As $\pi_2(K; \ast) = 0$, the map $g'$ gives rise to a cellular map of pairs $g: (X, \partial X) \longrightarrow (K,\partial Y)$ 
with $g |_{\partial X} = \tilde g$. Let
$g_{\ast}: C (X, \partial X) \longrightarrow C(K,\partial Y)$ be the chain map 
induced by $g$. Then $\varphi_{\ast}(c_{X \ast} ([X, \partial X])) = c_{Y \ast}([Y, \partial Y])$ 
implies $1 {\otimes}g_{\ast}(e) - 1 {\otimes}e' = 0$ in ${\mbox{H}}_{3}(K,\partial Y; {\mathbb{Z}}^{\omega_{Y}})$, where we 
identify $e'$ with its inclusion in $K$. Thus $g_{\ast}(e) - e'$ is contained in ${\rm{im}} (\partial_{3}: 
C_{4}(K, \partial Y) \rightarrow C_{3}(K, \partial Y)) + \overline{I_{Y}} C_{3}(K, \partial Y)$, and hence
\begin{equation}\label{IbarboundC3}
\partial_{2}(g_{\ast}(e) - e') \in \overline{I_{Y}} \partial_{2}(C_{3}(K, \partial Y)).
\end{equation}
The group isomorphism $\varphi: G \rightarrow H$ induces a ring
isomorphism $\varphi: {{\Lambda}}_{X} \rightarrow {{\Lambda}}_{Y}$ and $\varphi^{\ast}(\omega_{Y}) =
\omega_{X}$ implies $\varphi(\overline{{{\lambda}}}) =
\overline{\varphi({{\lambda}})}$ for every ${{\lambda}} \in {{\Lambda}}_{X}$. By Lemma~\ref{splitoffandmore},  $\overline{I_{X}}$ is generated by $\{ a_{m} \}_{m \in M}$, where $\partial_{2}(e) =
\sum_{m \in M} a_{m} e_{m}$. Hence $\overline{I_{Y}} = \overline{(\varphi(I_{X}))} = \varphi (\overline{I_{X}})$ is generated as right ${{\Lambda}}$--module by $\{ \varphi(a_{m}) \}_{m \in M}$
and we obtain
\[
\partial_{2}(g_{\ast}(e)) - \partial_{2}(e') =
g_{\ast}(\partial_{2}(e)) - \partial_{2}(e') = \sum_{m \in M}
\sum_k \varphi(a_{m}) {{\lambda}}_{mk} d_k,
\]
where $d_k = \partial_2 c_{k}$ for some $c_{k} \in C_{3}(K, \partial Y)$. In other words,
$d_{k}$ is represented by the gluing map $f_{k}: S^{2} \rightarrow 
q^{-1}(Y'\backslash \partial Y)$ of a $3$--cell in $q^{-1}(Y'\backslash \partial Y)$. 
We obtain $\partial_{2}(e') = \sum_{m} \varphi(a_{m}) (g_{\ast}(e_{m}) - \sum {{\lambda}}_{mk} d_{k})$
and as $D^{2} \# S^{2} \simeq D^{2}$, we may modify
$g'|_{X^{[2]}}$ on the interior of two--cells to obtain a map of pairs
$h: (X',\partial X) \longrightarrow (Y',\partial Y)$ with
\begin{equation}\label{threecellhom}
h_{\ast}(\partial_{2}(e)) = h_{\ast}(\sum_{m \in M} a_{m} e_{m})
= \sum_{m \in M} \varphi(a_{m}) (g_{\ast}(e_{m}) - \sum {{\lambda}}_{km} d_{k}) = \partial_{2}(e').
\end{equation}
Since all components of $\partial Y$ are aspherical, the components of
$q^{-1}(\partial Y)$ are open surfaces by Lemma~\ref{asphericalcompopensurf}.
Hence the long exact homology sequence of the pair $(Y', \partial Y)$
yields an injective homomorphism ${\mbox{H}}_{2}(q^{-1}(Y'); {{\Lambda}}_{Y}) \rightarrowtail 
{\mbox{H}}_{2}(q^{-1}(Y'), q^{-1}(\partial Y); {{\Lambda}}_{Y})$. As $\delta_{\ast} [Y,\partial Y] = [\partial Y]$, the class $[\partial_{2} e']$ is contained in the image of this
homomorphism. Therefore $\pi_{2}(Y') \cong {\mbox{H}}_{2}(q^{-1}(Y'); {{\Lambda}}_{Y})$
and (\ref{threecellhom}) imply that the composition of
the attaching map of the $3$--cell $e^{3}$ (represented by $e$) and the map
$h$ is homotopy equivalent to the attaching map of the $3$--cell
${e'}^{3}$ (represented by $e'$).

We conclude that $h$ extends to a map of pairs $f: (X,\partial X) \longrightarrow (Y,\partial Y)$
satisfying (i) -- (iii), showing that the ${\rm{PD}}^{3}$--pairs $(X,\partial X)$ and
$(Y, \partial Y)$ are homotopy equivalent. 
\end{proof}

\section{Realisation of Invariants by ${\rm{PD}}^3$--Pairs}
We need the following result from the \emph{derived module category}, also called the 
\emph{projective homotopy category of modules}.
\begin{theorem}\label{hefactorsas}
Let ${{\Lambda}}$ be a ring with unit. A homotopy equivalence $\varphi: A \rightarrow B$ of ${{\Lambda}}$--modules factors as
\[
\xymatrix{
A \ \ar@{>->}[r]^-{\iota} & A \oplus P \ar@{->>}[r]^-{\tilde{\varphi}}
& B \oplus Q \ar@{->>}[r]^-{\pi} & B,}
\]
where $P$ and $Q$ are projective and $\iota$ and $\pi$ are the
natural inclusion and projection respectively.
\end{theorem}
\begin{proof}
The proof is dual to the proof of Theorem 13.7 in \cite{Hilton}. 
\end{proof}
\begin{observation}
There is a projective ${{\Lambda}}$--module $\tilde P$ such that $P \oplus \tilde P$ is isomorphic to a free ${{\Lambda}}$--module $F$, so that  $\varphi$ factors as
\[
\xymatrix{
A \ \ar@{>->}[r] & A \oplus F \ar@{->>}[r] & B \oplus \tilde Q \ar@{->>}[r] & B }
\]
where $\tilde Q = Q \oplus \tilde P$ is projective. If the ${{\Lambda}}$--modules $A$ and 
$B$ are finitely generated, the projective ${{\Lambda}}$--modules $P$ and $Q$ are 
also finitely generated. Then $F \cong {{\Lambda}}^{n}$ for some $n \in {\mathbb{N}}$ and $\tilde Q$ is 
finitely generated projective.
\end{observation}
Let ${{\Lambda}}$ be the integral group ring of the group $H$ and take $\omega \in {\mbox{H}}^{1}(H; {\mathbb{Z}} / 2 {\mathbb{Z}})$.
Given a chain complex $\ldots \rightarrow C_{r+1} \stackrel{\partial_{r}}{\rightarrow} C_{r}
\rightarrow \ldots$ of left ${{\Lambda}}$--modules, put
\[
{\rm{G}}_{r}(C) := {\rm{coker}}\partial_{r} = C_{r} / {\rm{im}}\partial_{r},
\]
and given a chain map $f: C \rightarrow D$, let ${\rm{G}}_{r}(f): {\rm{G}}_{r}(C) \rightarrow
{\rm{G}}_{r}(D)$ be the induced ${{\Lambda}}$--morphism of cokernels. Then $G = G_{\ast}$ is a functor 
from the category of chain complexes of left ${{\Lambda}}$--modules to itself. Writing $C^{\ast}$ for $^{\omega}\!{{\rm{Hom}}_{\Lambda}(C, {{\Lambda}})}$, we compose the functors $G$ and $^{\omega}\!{{\rm{Hom}}_{\Lambda}(-, {{\Lambda}})}$ 
to obtain the functor $F$ (see \cite{Turaev} p.265) given by
\begin{equation}\label{Fdfn}
F^{r}(C) = G_{-r}(C^{\ast}) = C^{r}/{\rm{im}}\partial_{r-1}^{\ast}.
\end{equation}
\begin{lemma}\label{chaintostab}
Let $f,g: C \rightarrow D$ be chain homotopic maps of chain complexes
over ${{\Lambda}}$. If $D_{n}$ is projective, then $G_{n}(f) \simeq G_{n}(g)$ as
${{\Lambda}}$--morphisms.
\end{lemma}
For a proof we refer the reader to \cite{Bleile}.

Thus we may view $G_{n}$ as a functor from the category of chain complexes of projective left
${{\Lambda}}$--modules and homotopy classes of chain maps to the derived module category.
\begin{corollary}\label{Fhomeq}
Let $f: C \rightarrow D$ be a homotopy equivalence of chain
complexes over ${{\Lambda}}$. If $C_{n}$ and $D_{n}$ are projective, then
$G_{n}(f)$ is a homotopy equivalence of ${{\Lambda}}$--modules.
\end{corollary}
\begin{lemma}\label{compcapwithgen}
    Let $(X,Y)$ be a pair of $CW$--complexes with $X$
    connected, $\omega \in {\mbox{H}}^{1}(X;{\mathbb{Z}} / 2 {\mathbb{Z}})$ and
    ${\mbox{H}}_{n}(X,Y;{\mathbb{Z}}^{\omega}) \cong {\mathbb{Z}}$ generated by $[1 {\otimes}x]$.
    Then there is a chain $w_{1} \in C_{1}(X)$, such that the
    ${{\Lambda}}$--morphism $\cap 1 {\otimes}x: C^{\ast}(X,Y)
    \rightarrow ^{\omega}\!\!{{{\Lambda}}}^{\omega} {\otimes_{\Lambda}} C(X) \cong C(X)$ is given by
    \[ \varphi \cap 1 {\otimes}x = \overline{\varphi(x)}.(1 + \partial_{0} w_{1}) \]
    for every cocycle $\varphi \in C^{\ast}(X,Y)$, where we
    identify ${{\lambda}} {\otimes}c \in \, ^{\omega}\!{{{\Lambda}}}^{\omega} {\otimes_{\Lambda}} C(X)$ with
    $\overline{{{\lambda}}}.c \in C(X)$.
\end{lemma}
\begin{proof}
    Let $\pi: C(X) \twoheadrightarrow C(X, \partial X)$ be the natural projection,
    take $y \in C_{n}(X)$ with $\pi(y) = x$ and assume $\Delta y = \sum
    y_{i} {\otimes}z_{n-i}$ with $y_{i}, z_{i} \in C_{i}(X)$. Then $y = ({\rm{id}}
    {\otimes}\varepsilon)\Delta(y)$ implies $y_{n}.\varepsilon(z_{0}) =
    y$. As $[1 {\otimes}x]$ is a generator, $x$ and thus $y$ are
    indivisible, so that $y = y_{n}$ and $\varepsilon(z_{0}) = 1$ up to sign. As
    $X$ is connected, we may assume $C_{0}(X) = {{\Lambda}}$ and identify
    ${\rm{im}}\partial_{1}$ with $I = \ker \varepsilon$. Then
    $\varepsilon(z_{0}) = 1$ implies $z_{0} = 1 + w_{0}$ with $w_{0}
    \in I$. Hence $z_{0} = 1 + \partial_{0} w_{1}$ for some
    $w_{1} \in C_{1}(X)$, so that $\varphi \cap 1 {\otimes}x =  
    \overline{\varphi(x)}.(1 + \partial_{0} w_{1}).$
\end{proof}

Given a ${\rm{PD}}^{3}$--pair, $(X, \partial X)$, take a cycle 
$1 {\otimes}x \in {\mathbb{Z}}^{\omega} {\otimes_{\Lambda}} C_{3}(X, \partial X)$ representing $[X, \partial X]$. Then
$\cap 1 {\otimes}x: C^{\ast}(X,\partial X)\rightarrow^{\omega}\!\!{{{\Lambda}}}^{\omega} {\otimes_{\Lambda}} C(X) \cong C(X)$
is a chain homotopy equivalence. As both $C_{2}^{\ast}(X,\partial X)$ and $C_{1}(X)$ are free 
and hence projective, Corollary~\ref{Fhomeq} implies that
\[ G_{-2}(\cap 1 {\otimes}x): F^{2}(C(X, \partial X)) = G_{-2}(C^{\ast}(X,\partial X))
\rightarrow G_{1}(C(X))\]
is a homotopy equivalence of ${{\Lambda}}$--modules. Composing with the isomorphism
$\vartheta: G_{1}(C(X)) \rightarrow I$, given by $\vartheta([c]) := \partial_{0}(c)$,
we obtain another homotopy equivalence of ${{\Lambda}}$--modules. The fact that 
$\cap 1 {\otimes}x$ is a chain map together with Lemma~\ref{compcapwithgen} yields
$(\vartheta \circ G_{-2}(\cap 1 {\otimes}x) ) ([\varphi]) = \overline{(\partial_{2}^{\ast} \varphi)(x)} 
+ \overline{(\partial_{2}^{\ast} \varphi)(x)}\partial_{0} w_{1}$
for every $\varphi \in C_{2}^{\ast}(X, \partial X)$ and some $w_{1} \in C_{1}(X)$. As the ${{\Lambda}}$--morphism $F^{2}(C(X, \partial X)) \longrightarrow I, \ [\varphi] \longmapsto
\overline{(\partial_{2}^{\ast} \varphi)(x)}\partial_{0} w_{1}$ is null--homotopic,
\[F^{2}(C(X, \partial X)) \longrightarrow I, 
\ [\varphi] \longmapsto \overline{(\partial_{2}^{\ast} \varphi)(x)}\]
is a homotopy equivalence of ${{\Lambda}}$--modules.

Attaching cells of dimension three and larger to $(X, \partial X)$,
we obtain an Eilenberg--Mac\,Lane pair $(K, \partial X)$ of
type $K( \{\kappa_{i}: \pi_{1}(\partial X_{i}, \ast) \rightarrow \pi_{1}
(X, \ast) \}_{i \in J};1)$ with the inclusion $\iota: (X, \partial X)
\rightarrow (K, \partial X)$ a cellular classifying map. Identifying the 
cellular chain complexes of the pair $(X, \partial X)$ with their image under the 
chain map induced by $\iota$, we obtain $C_{i}(K) = 
C_{i}(X), C_{i}(K, \partial X) = C_{i}(X, \partial X)$ for $i = 0, 1, 2$ and 
$[1 {\otimes}x] = [X, \partial X] = \iota_{\ast}([X, \partial X])$. 
\begin{lemma}\label{itsaneq}
    The ${{\Lambda}}$--morphism
    \begin{equation}\label{crucialhomeq2}
        F^{2}(C(K, \partial X)) \longrightarrow I, \ [\varphi] \longmapsto
        \overline{(\partial_{2}^{\ast} \varphi)(x)}
    \end{equation}
    is a homotopy equivalence of ${{\Lambda}}$--modules.

\end{lemma}
Given a chain complex $C$ of free left ${{\Lambda}}$--modules, Turaev
constructed a group homomorphism
\[ \nu_{C,r}: {\mbox{H}}_{r+1}({\mathbb{Z}}^{\omega} {\otimes_{\Lambda}} C) \longrightarrow [F^{r},I] \]
such that $\nu_{C(X, \partial X),2}([1 {\otimes}x]) = \nu_{C(K, \partial
X),2}(\iota_{\ast}([X, \partial X]))$ is the homotopy class of the
homotopy equivalence~(\ref{crucialhomeq2}). Let $\delta$ denote the 
connecting homomorphism of the short exact sequence 
$\overline{I} C \rightarrowtail C \twoheadrightarrow {\mathbb{Z}}^{\omega} {\otimes_{\Lambda}} C$
and identify $c \in C_{r}$ with $1 {\otimes}c \in {{\Lambda}}^{\omega} {\otimes_{\Lambda}} C$. Then the
natural transformation $\eta$ of (\ref{eta}) yields the ${{\Lambda}}$--morphism
$\eta: C_{r} \longrightarrow (C_{r}^{\ast})^{\ast}, \, c \longmapsto \eta(c)$
given by $\eta(c) (\varphi) = \overline{\varphi(c)}$. It is not difficult to show that
$\eta(c)$ factors through the cokernel $F^{r}(C)$ of $\partial^{\ast}_{r-1}$ and
that its image is contained in $I$. We obtain the ${{\Lambda}}$--morphism
\[ \hat{\eta(c)}: F^{r}(C) \longrightarrow I, \, [\varphi]
\longmapsto \overline{\varphi(c)},\]
whose homotopy class depends on the homology
class of the cycle $c \in \overline{I} C$ only. The composition of the homomorphism
${\mbox{H}}(\overline{I} C) \longrightarrow [F^{r}(C),I], \, [c] \longmapsto [\hat{\eta(c)}]$ 
with $\delta$ yields the homomorphism
\begin{equation}\label{nudef}
    \nu_{C,r}: {\mbox{H}}_{r+1}({\mathbb{Z}}^{\omega} {\otimes_{\Lambda}} C) \longrightarrow [F^{r}(C),I], \quad
    [1 {\otimes}c] \longmapsto [ \hat{\eta(c)}]
\end{equation}
represented by the ${{\Lambda}}$--morphism
  \[ F^{r}(C) \longrightarrow I, \, [\varphi] \longmapsto
    \overline{\varphi(\partial_{r}(x))}.\]
\begin{lemma}\label{nuisinj}
    Suppose $C$ is a chain complex of free left ${{\Lambda}}$--modules
    such that $C_{r}$ is finitely generated and ${\mbox{H}}_{r}(C) =
    {\mbox{H}}_{r+1}(C) = 0$. Then $\nu_{C,r}$ is an isomorphism.
\end{lemma}
For a proof we refer the reader to Lemma 2.5 in \cite{Turaev}.

Suppose $(K,\partial X)$ is an Eilenberg--Mac\,Lane pair of type 
$K(\{\kappa_{i}: G_i \rightarrow G\}_{i \in J};1)$, 
$\omega \in {\mbox{H}}^{1}(K, {\mathbb{Z}} / 2 {\mathbb{Z}})$ and 
$\mu \in \mbox{H}_3(K, \partial X; {\mathbb{Z}}^{\omega})$.
\begin{theorem}[REALSIATION I]\label{neccon}
    If $(\{\kappa_{i}\}_{i \in J}, \omega, \mu)$ is
    the fundamental triple of a ${\rm{PD}}^{3}$--pair
    then 
    $\nu_{C(K, \partial X),2}(\mu)$ is a homotopy 
    equivalence of ${{\Lambda}}$--modules.
\end{theorem}
\begin{proof}
    Assume $(\{\kappa_{i}\}_{i \in J}, \omega, \mu)$ is
    the fundamental triple of the ${\rm{PD}}^{3}$--pair $(X, \partial X)$ and the
    Eilenberg--Mac\,Lane pair $(K,\partial X)$ of type $K(\{\kappa_{i}\}_{i \in J}; 1)$ 
    was obtained by attaching cells of dimension three and larger to $X$.
    Take 
    \[ 1 {\otimes}x \in {\mathbb{Z}}^{\omega} {\otimes_{\Lambda}} C_{3}(X,
    \partial X) \subseteq {\mathbb{Z}}^{\omega} {\otimes_{\Lambda}} C_{3}(K, \partial X)\]
    with $[1 {\otimes}x] = \mu$.
    Then $F^{2}(C(K, \partial X)) \longrightarrow I, \, [\varphi]
    \longmapsto \overline{\varphi(\partial_{2}(x))}$ represents the class $\nu_{C(K, \partial X),2}(\mu)$ and is a homotopy equivalence of ${{\Lambda}}$--modules    by  
    Lemma~\ref{itsaneq}.

Given another Eilenberg--Mac\,Lane pair $(L, \partial L)$ of type $K(\{\kappa_{i}\}_{i \in J}; 1)$,
there is a homotopy equivalence $f: (K, \partial X) \rightarrow (L, \partial L)$ of pairs of 
CW--complexes inducing a chain homotopy equivalence 
$g^{\ast}: C^{\ast}(K, \partial X) \longrightarrow C^{\ast} (L, \partial L)$. Thus Corollary~\ref{Fhomeq}
implies that $F^{2}(g) = G_{-2}(g^{\ast})$ is a homotopy equivalence of ${{\Lambda}}$--modules. The diagram
    \[ \xymatrix{
    F^{2}(C(L, \partial L)) \ar[d]_{F^{2}(g)} \ar[rrr]^{\nu_{C(L, \partial
    L),2}(f_{\ast}\mu)} &&& I \ar@{=}[d] \\
    F^{2}(C(K, \partial K)) \ar[rrr]^{\nu_{C(K, \partial
    K),2}(\mu)} &&& I }\]
    commutes and hence $\nu_{C(L, \partial L),2}(f_{\ast}\mu)$ is a
    homotopy equivalence of ${{\Lambda}}$--modules if and only if $\nu_{C(K, \partial
    K),2}(\mu)$ is one.
\end{proof}

For  $\{\kappa_{i}: G_{i} \rightarrow G\}_{i \in J}$ to be the
$\pi_{1}$--system of a ${\rm{PD}}^{3}$--pair $(X, \partial X)$, the groups
$G_{i}$ must be surface groups for all $i \in J$ as the components
of $\partial X$ are ${\rm{PD}}^{2}$--complexes by definition and thus homotopy
equivalent to closed surfaces.  

\begin{definition}
The $\pi_1$--system $\{\kappa_{i}\}_{i \in J}$ is \emph{injective} if $\kappa_{i}$ is injective 
for every $ i \in J$. The triple $(\{\kappa_{i}\}_{i \in J},\omega, \mu)$ is then called 
\emph{$\pi_1$--injective}. 
\end{definition}
Let $\{\kappa_{i}: G_{i} \rightarrow G\}_{i \in J}$ be an injective $\pi_{1}$--system
with $G$ finitely presentable and $G_{i}$ a surface group for all $i \in J$. Furthermore, let 
$(K, \partial X)$ be an Eilenberg--Mac\,Lane pair of type $K(\{\kappa_{i}\}_{i \in J};1)$
such that the components $\partial X_{i}$ of $\partial X$ are all surfaces.
Take $\omega \in {\mbox{H}}^{1}(K; {\mathbb{Z}} / 2 {\mathbb{Z}})$ and 
$\mu \in {\mbox{H}}_{3}(K, \partial X; {\mathbb{Z}}^{\omega})$ such that 
$\delta_{\ast} \mu = [\partial X]$, where $[\partial X]$ is the fundamental class 
of the ${\rm{PD}}^{2}$--complex $\partial X$ and $\delta_{\ast}$ is the connecting 
homomorphism of $C(\partial X) \rightarrowtail C(X) \twoheadrightarrow C(X, \partial X)$.
Following Turaev's construction in the absolute case, we construct a
${\rm{PD}}^{3}$--pair realising $(\{\kappa_{i}\}_{i \in J}, \omega, \mu)$.

Since $G$ is assumed finitely presentable, we may also assume that $K$ has 
finite 2--skeleton $K^{[2]}$, so that the ${{\Lambda}}$--modules $C_{2}(K, \partial X)$ 
and thus $F^{2}(C(K, \partial X))$ are finitely generated.  Let 
$h: F^{2}(C(K, \partial X)) \rightarrow I$ be a ${{\Lambda}}$--morphism representing 
$\nu_{C(K,\partial X),2}(\mu)$. Then $h$ is a homotopy equivalence of ${{\Lambda}}$--
modules which factors as
\[ \xymatrix{
F^{2}(C(K, \partial X)) \,\, \ar@{>->}[r] &
F^{2}(C(K, \partial X)) \oplus {{\Lambda}}^m \,\, \ar@{>->>}[r]^-{j} &
I \oplus P \ar@{->>}[r] & I},\]
by Theorem~\ref{hefactorsas}, where $P$ is finitely generated and
projective. Let $B = (e^{0} \vee e^{2}) \cup e^{3}$ be a $3$--dimensional ball. If
we replace $K$ by $K \vee (\bigvee_{i = 1}^m B)$, then $K^{[2]}$ is
replaced by $K^{[2]} \vee (\bigvee_{i = 1}^m e^{2})$ and $F^{2}(C(K, \partial
X))$ is replaced by $F^{2}(C(K, \partial X)) \oplus {{\Lambda}}^m$. 
Thus we may assume without loss of generality that $h$ factors as
\begin{equation}\label{hfac} 
\xymatrix{
F^{2}(C(K, \partial X)) \,\, \ar@{>->>}[r]^-{j} &
I \oplus P \ar@{->>}[r] & I,}
\end{equation}
where $P$ is finitely generated and projective.

First we consider the case where $P$ is free, that is, $P \cong {{\Lambda}}^n$, 
for some $n \in {\mathbb{N}}$. Denote the natural projection $C^{2}(K, \partial X) 
\twoheadrightarrow F^{2}(C(K, \partial X))$ by $\pi$, put $\tilde \iota := (\iota, {\rm{id}}_P): 
I \oplus P \rightarrow {{\Lambda}} \oplus P$ and  consider the ${{\Lambda}}$--morphism
\begin{equation}
    \xymatrix{
    \varphi: C^{2}(K, \partial X) \ar@{->>}[r]^-{\pi} &
    F^{2}(C(K, \partial X)) \,\, \ar@{>->>}[r]^-{j} &
    I \oplus P \,\, \ar@{>->}[r]^-{\tilde \iota} & {{\Lambda}} \oplus P.}
\end{equation}
By definition, $\varphi \circ \partial_{1}^{\ast} = 0$ and hence 
${\rm{im}} \varphi^{\ast} \subseteq \ker \partial_{1}$.

Let $p: \tilde K \rightarrow K$ be the universal covering. Since
$\kappa_{i}$ is injective for every $i \in J$, the components of
$p^{-1}(\partial X)$ are universal covering spaces of
Eilenberg--Mac\,Lane complexes, so that ${\mbox{H}}_{2}(p^{-1}(\partial X)) =
{\mbox{H}}_{1}(p^{-1}(\partial X)) = 0$. The long exact homology sequence
of the pair $(p^{-1}(K^{[2]}), p^{-1}(\partial X))$ and the Hurewicz 
Isomorphism Theorem imply
\begin{equation}\label{needinj}
{\rm{im}} \varphi^{\ast} \subseteq \ker \partial_{1}
=         {\mbox{H}}_{2}(p^{-1}(K^{[2]}), p^{-1}(\partial X)) \\
\cong      {\mbox{H}}_{2}(p^{-1}(K^{[2]})) \\
\cong      \pi_{2}(p^{-1}(K^{[2]})).
\end{equation}
We may thus attach $3$--cells to $K^{[2]}$ to
obtain a pair $(X, \partial X)$ of CW--complexes whose relative
cellular chain complex is given by
\[ D: 0 \longrightarrow ({{\Lambda}} \oplus P)^{\ast} \stackrel{\varphi^{\ast}}
{\longrightarrow} C_{2}(K, \partial X) \longrightarrow
C_{1}(K, \partial X) \longrightarrow C_{0}(K, \partial X).\]
As $\pi_{2}(K) = 0$, the inclusion 
\[ (K^{[2]}, \partial X) \rightarrow (K, \partial X)\]
extends to a map $f: (X, \partial X) 
\longrightarrow (K, \partial X)$. Since $f$ induces an isomorphism of
$\pi_{1}$--systems, we may view $\omega$ as an element of
${\mbox{H}}^{1}(X; {\mathbb{Z}} / 2 {\mathbb{Z}})$.

Then $(X, \partial X)$ is a ${\rm{PD}}^{3}$--pair realising 
$(\{ \kappa_{i} \}_{i \in J}, \omega, \mu)$ since
\begin{proposition}
\hfill
    \begin{itemize}
    \item[(i)] ${\mbox{H}}_{3}(X, \partial X; {\mathbb{Z}}^{\omega}) \cong {\mathbb{Z}}$;
    \item[(ii)] $f_{\ast}([X, \partial X]) = \mu$, where
                    $[X, \partial X]$ generates ${\mbox{H}}_{3}(X, \partial X; {\mathbb{Z}}^{\omega})$;
    \item[(iii)] $\delta_{\ast} [X, \partial X] = [\partial X]$, where
                     $[\partial X]$ is the fundamental class of the
                     ${\rm{PD}}^{2}$--complex $\partial X$ and $\delta_{\ast}$ is the
                     connecting homomorphism of the short exact sequence
                     \[ C(\partial X) \rightarrowtail C(X) \twoheadrightarrow
                     C(X, \partial X);\]
     \item[(iv)] $\cap [X, \partial X]: {\mbox{H}}^{r}(X; ^{\omega}\!{{{\Lambda}}}^{\omega})
                      \rightarrow {\mbox{H}}_{r-3}(X, \partial X; {{\Lambda}})$ is an
                      isomorphism for every $r \in {\mathbb{Z}}$.
     \end{itemize}
\end{proposition}
\begin{proof}
We only prove (iv), details for (i) -- (iii) are contained in \cite{Bleile}.

First observe that the definition of $(X, \partial X)$ implies
${\mbox{H}}^{2}(X, \partial X; ^{\omega}\!{{{\Lambda}}}^{\omega}) = 0.$ Since 
${\mbox{H}}_{1}(X; {{\Lambda}}) = {\mbox{H}}_{1}(C(X)) = 0$, the homomorphism
\[ \cap [X, \partial X]: {\mbox{H}}^{2}(X, \partial X; ^{\omega}\!{{{\Lambda}}}^{\omega}) \rightarrow
{\mbox{H}}_{1}(X; {{\Lambda}})\]
is an isomorphism.

As ${{\Lambda}} {\otimes}P$ is free, we may use the natural
transformation $\eta$ to identify the module
$^{\omega}\!{\rm{Hom}}_{\Lambda}(({{\Lambda}} \oplus P)^{\ast}, ^{\omega}\!{{{\Lambda}}}^{\omega})$ 
with ${{\Lambda}} \oplus P$ and $(\varphi^{\ast})^{\ast}$ with $\varphi$. Then
${\mbox{H}}^{3}(X, \partial X; ^{\omega}\!{{{\Lambda}}}^{\omega}) \cong  ({{\Lambda}} \oplus P) / {\rm{im}} \varphi  
\cong  {{\Lambda}} / I \cong  {\mathbb{Z}}.$ Clearly, ${\mbox{H}}^{3}(X, \partial X; ^{\omega}\!{{{\Lambda}}}^{\omega})$ is generated by
$\psi = (1, 0) \in ({{\Lambda}}^{\ast})^{\ast} \oplus (P^{\ast})^{\ast} =
C_{3}^{\ast}(X, \partial X) = C_{3}^{\ast}(X)$. By Lemma~\ref{compcapwithgen},
$[\psi] \cap [X, \partial X] = [\psi] \cap [1 {\otimes}x] = \overline{\psi(x)} = 1$, that is, 
$\cap [X, \partial X]$ maps $\psi$ to a generator of ${\mbox{H}}_{0}(X; {{\Lambda}})$. Hence
$\cap [X, \partial X]: {\mbox{H}}^{3}(X, \partial X; ^{\omega}\!{{{\Lambda}}}^{\omega}) \rightarrow
{\mbox{H}}_{0}(X; {{\Lambda}})$ is an isomorphism. 

Since $\partial X$ is a ${\rm{PD}}^{2}$--complex, $\cap [\partial X]: {\mbox{H}}^{r}(\partial X; ^{\omega}\!{{\Lambda}}^{\omega}) \longrightarrow {\mbox{H}}_{2-r}(\partial X; {{\Lambda}})$ is an isomorphism 
for every $r \in {\mathbb{Z}}$. Thus the Cap Product Ladder of $(X, \partial X)$ (Theorem \ref{ladder})
with $y = [X, \partial X]$, and the Five Lemma imply that
$\cap [X, \partial X]: {\mbox{H}}^{r}(X; ^{\omega}\!\!{{{\Lambda}}}^{\omega}) \rightarrow
{\mbox{H}}_{r-3}(X, \partial X; {{\Lambda}}) $ is an isomorphism for $r = 2$ and $r =3$. Identifying $^{\omega}\!\!{{{\Lambda}}}^{\omega}$ with ${{\Lambda}}$, ${{\Lambda}} {\otimes_{\Lambda}} A$ with $A$ and denoting 
$^{\omega}{\rm{Hom}}_{\Lambda}(A, {{\Lambda}})$ by $A^{\ast}$ for left ${{\Lambda}}$--modules $A$, we obtain the 
chain homotopy equivalence
\begin{equation*}
    \xymatrix{
    0 \ar[r] \ar[d] & {\rm{im}} \partial_{1}^{\ast} \ \ar@{>->}[r] \ar[d]&
    C_{2}^{\ast}(X) \ar[r] \ar[d]^{\cap 1 {\otimes}x} &
    C_{3}^{\ast}(X) \ar[r] \ar[d]^{\cap 1 {\otimes}x} & 0 \ar[d] \\
    0 \ar[r] & {\rm{im}} \partial_{1} \ \ar@{>->}[r]  &
    C_{1}(X, \partial X) \ar[r]  &
    C_{0}(X, \partial X) \ar[r]  & 0.}
\end{equation*}
Applying $^{\ast} = \ ^{\omega}{\rm{Hom}}_{\Lambda}(-, {{\Lambda}})$
yields the chain homotopy equivalence
\begin{equation*}
\xymatrix{
0 \ar[r] \ar[d] & C_{0}(X, \partial X)^{\ast} \ar[r] \ar[d]^{(\cap 1 {\otimes}x)^{\ast}} & 
C_{1}(X, \partial X)^{\ast} \ar[r] \ar[d]^{(\cap 1 {\otimes}x)^{\ast}}  &
({\rm{im}} \partial_{1})^{\ast} \ar[r] \ar[d]  & 0 \ar[d]  \\
0 \ar[r] & (C_{3}(X)^{\ast})^{\ast} \ar[r] & (C_{2}(X)^{\ast})^{\ast} \ar[r] &
({\rm{im}} \partial_{1}^{\ast})^{\ast} \ar[r] & 0,}
\end{equation*}
which shows that $(\cap [1 {\otimes}x])^{\ast}$ induces homology isomorphisms.
By Lemma~\ref{capastcom}, the $\Lambda$--morphism $\cap(1 {\otimes}x)$ induces isomorphisms in 
homology if and only if $(\cap 1 {\otimes}x)^{\ast}$ does. Thus the homomorphism
$\cap [X, \partial X]: {\mbox{H}}^{k}(X, \partial X;
^{\omega}\!\!{{{\Lambda}}}^{\omega}) \longrightarrow {\mbox{H}}_{3-k}(X; {{\Lambda}})$ is an isomorphism for 
$k =0$ and $k =1$. The Cap Product Ladder of $(X, \partial X)$ with $y = [X,
\partial X]$ and the Five Lemma imply that
$\cap [X, \partial X] : {\mbox{H}}^{r}(X; ^{\omega}\!\!{{{\Lambda}}}^{\omega}) \longrightarrow
{\mbox{H}}_{3-k}(X, \partial X; {{\Lambda}})$ is an isomorphism for $r =0$ and $r = 1$ and 
hence for every $r \in {\mathbb{Z}}$.
\end{proof}
\begin{theorem}[REALISATION II]
There is a ${\rm{PD}}^{3}$--pair which realises the $\pi_1$--injective triple $(\{\kappa_{i}\}_{i \in J},\omega, \mu)$ if and only if $\nu_{C(K,\partial X),2}(\mu)$ is a class of homotopy eqivalences.
\end{theorem}
\begin{proof}
It only remains to investigate the general case, where the module $P$ in the 
factorisation~(\ref{hfac}) of the homotopy equivalence $h$ is 
finitely generated projective, but not necessarily free. Then there is 
a finitely generated projective ${{\Lambda}}$--module $Q$ such that $P^{\ast} \oplus Q = {{\Lambda}}^{n}$. Attaching infinitely many $3$--cells to 
$K^{[2]} \vee \big( \vee_{i=1}^{\infty} e^{2} \big)$ we obtain a pair 
$(X, \partial X)$ of $CW$--complexes whose relative cellular chain 
complex is chain homotopy equivalent to the complex
\[
\xymatrix{
E: \quad \ldots \ar[r] & 
{{\Lambda}}^{n} \ar[r]^{{\rm{pr}}} & 
{{\Lambda}}^{n} \ar[r]^{{\rm{pr'}}} & 
{{\Lambda}}^{n} \ar[r]^-{q}  &
({{\Lambda}} \oplus P)^{\ast} \oplus Q \\
\ar[rr]^-{\begin{bmatrix} \varphi^{\ast} & 0 \\ 0  & 0 \end{bmatrix}} &&
C_{2}(K, \partial X) \ar[r] & 
C_{1}(K, \partial X) \ar[r] & 
C_{0}(K, \partial X), & &}\]
where ${\rm{pr}}: {{\Lambda}}^{n} = P^{\ast} \oplus Q \rightarrow Q$ and
${\rm{pr'}}: {{\Lambda}}^{n} = P^{\ast} \oplus Q \rightarrow P^{\ast}$ are the 
canonical projections and $q(x) = (0,0,{\rm{pr}}(x)) \in 
({{\Lambda}} \oplus P)^{\ast} \oplus Q$ for $x \in {{\Lambda}}^{n}$. As $E$ is a complex of 
finitely generated free ${{\Lambda}}$--modules, the proof that $(X, \partial X)$ 
realises $(\{ \kappa_{i} \}_{i \in J}, \omega, \mu)$ is analogous to the 
proof in the case where the module $P$ is free. 
\end{proof}
\begin{observation}
If the ${\rm{PD}}^{3}$--pair $(X, \partial X)$ realises $(\{\kappa_{i}\}_{i \in J},\omega, \mu)$ with $G$ finitely presentable, Wall's results imply that $X$ is in fact dominated by a finite cell--complex.
\end{observation}
\begin{observation}
Note that  $\pi_1$--injectivity guarantees the first isomorphism in~(\ref{needinj}) which allows us to attach $3$--cells to $K^{[2]}$ such that $\varphi^{\ast}$ is the boundary operator of the resulting relative cellular chain complex. Thus $\pi_1$--injectivity is an indispensable assumption for our method. The question remains whether the realisation theorem holds without this assumption.
\end{observation}
\section{Connected Sums and Decomposition of ${\rm{PD}}^3$--Pairs}
\subsection{The Interior Connected Sum}
Given a ${\rm{PD}}^3$--pair $(X, \partial X)$, we may assume that the components of $\partial X$ are collared in $X$, the collar being a $3$--manifold with boundary. Lemma~\ref{splitoffandmore} implies that there is a homologically $2$--dimensional $CW$--complex $K$ and a map $f: S^2 \rightarrow K$ such that
\[ X \sim K  \cup_f e^3.\]
We may assume without loss of generality that $f(S^2) \cap \partial X = \emptyset$. If $e \in C_{3}(X, \partial X)$ is the chain corresponding to the $3$--cell $e^{3}$, then $1 {\otimes}e$ is a relative cycle representing the homology class $[X,\partial X]$. Wall \cite{Wall1} showed that the pair $(K, f)$ is unique up to homotopy and orientation for $\partial X = \emptyset$. The argument extends to the case of ${\rm{PD}}^3$--pairs with $\partial X$ not necessarily empty, allowing for the definition of the interior connected sum of  ${\rm{PD}}^3$--pairs. 
\begin{proposition}\label{uniquesplit}
The pair $(K, f)$ is unique up to homotopy, that is, if $X \sim K_1 \cup_{f_1} e^3 \sim K_2 \cup_{f_2} e^3$, there is a homotopy equivalence $g: K_1 \rightarrow K_2$ such that $g_{\ast}([f_1]) = [f_2]$, where $g_{\ast}$ denotes the isomorphism of homotopy groups induced by $h$.
\end{proposition}
For a proof we refer the reader to \cite{Bleile}.

To define the interior connected sum of two ${\rm{PD}}^3$--pairs $(X_1, \partial X_1)$ and $(X_2, \partial X_2)$, write $X_{\ell} = K_{\ell} \cup_{f_{\ell}} e_{\ell}^3$ for ${\ell} = 1, 2,$ and let $\iota_{\ell}: K_{\ell} \rightarrow K_1 \vee K_2$ be the inclusion of the first and second factor respectively. Then $\hat f_{\ell} := \iota_{\ell} \circ f_{\ell}: S^2 \longrightarrow K_1 \vee K_2$ determines an element of $\pi_2(K_1 \vee K_2; \ast)$, and we put $f_1 + f_2 := \hat f_1 + \hat f_2$.
\begin{lemdfn}\label{intsumdfn}
Up to oriented homotopy equivalence, the pair 
\[ (X, \partial X) := (X_1, \partial X_1) \sharp (X_2, \partial X_2) := \big( (K_1 \vee K_2) \cup_{f_1+f_2} e^3, \partial X_1 \cup \partial X_2) \big) \]
is a ${\rm{PD}}^3$--pair uniquely determined by $(X_i, \partial X_i), i = 1, 2$. It is called the {\emph{interior connected sum}} of $(X_1, \partial X_1)$ and $(X_2, \partial X_2)$.
\end{lemdfn}
\begin{proof}
Proposition~\ref{uniquesplit} implies uniqueness up to oriented homotopy equivalence. Theorem~\ref{sumthm} guarantees that $(X, \partial X)$ is indeed a ${\rm{PD}}^3$--pair. To see this, put
$G := \pi_1(X; \ast)$ and $G_k :=  \pi_1(X_k; \ast)$, $k = 1, 2$, so that $G = G_1 \ast G_2$. For $k = 1, 2$, let $\iota_k: G_k \rightarrow G$ be the canonical inclusion.  Regarding ${\mathbb{Z}}[G]$ as a right ${\mathbb{Z}}[G_k]$--module via $\iota_k$, we define the functor $L_k$ from the category of ${\mathbb{Z}}[G_k]$--modules to the category of ${\mathbb{Z}}[G]$--modules by 
\begin{equation}\label{functor}
L_k M := {\mathbb{Z}}[G] \otimes_{{\mathbb{Z}}[G_k]} M \quad {\rm{and}} \quad L_k \alpha := {\rm{id}} \otimes \alpha
\end{equation} 
for ${\mathbb{Z}}[G_k]$--modules $M$ and ${\mathbb{Z}}[G_k]$--morphisms $\alpha: M \rightarrow N$. Let $B$ denote the subcomplex of $C(X)$ containing the $3$--cell attached via $f_1 + f_2$ together with its boundary and denote the boundary by $\partial B$. Then $(L_k(C(K_k)), L_k(C(\partial X_k)) + \partial B)$ is a Poincar\'{e} chain pair for $k = 1, 2$ and $(\partial B, \emptyset)$ is also a Poincar\'{e} chain pair. Applying Theorem~\ref{sumthm} repeatedly, we see that $(C(X), C(\partial X)) = ( L_1(C(K_1)) + L_2(C(K_2)) + B, L_1(C(\partial X_1)) + L_2(C(\partial X_2)))$ is a Poincar\'{e} chain pair, showing that $(X, \partial X)$ is a ${\rm{PD}}^3$--pair.
\end{proof}
The notion of interior connected sums of ${\rm{PD}}^3$--pairs is consistent with that of interior connected sums of manifolds with boundary. In the case of empty boundaries, the notion of interior connected sums of ${\rm{PD}}^3$--pairs reduces to that of connected sums of ${\rm{PD}}^3$--complexes \cite{Wall1}.

\begin{definition}
The free product of two $\pi_1$--systems $\{ \kappa_{k \ell}: G_{k \ell} \rightarrow G_k \}_{\ell \in J_k}$, 
$k = 1, 2$, is the $\pi_1$--system 
\begin{equation*}
\{ \iota_k \circ \kappa_{k \ell}: G_{k \ell} \rightarrow G_1 \ast G_2 \}_{\ell \in J_k, k = 1, 2},
\end{equation*}
where $\iota_k: G_k \rightarrow G_1 \ast G_2$ denotes the inclusion of the factor in the free product of groups. 

A $\pi_1$--system isomorphic to such a free product is said to \emph{decompose as free product}.
\end{definition}
A decomposition of the ${\rm{PD}}^3$--pair $(X, \partial X)$ as interior connected sum yields a decomposition of its $\pi_1$--system as free product. The Classification Theorem \ref{B} and the Realisation Theorem \ref{realthm} allow us to show that the converse holds for finitely dominated pairs with non--empty aspherical boundaries in the case of $\pi_1$--injectivity, thus proving Theorem \ref{decomp1}.
\begin{proof}[Proof of Decomposition I]
Suppose the $\pi_1$--system of the finitely dominated ${\rm{PD}}^3$--pair $(X, \partial X)$ with aspherical boundary components decomposes as free product of the injective $\pi_1$--systems $\{ \kappa_{k \ell}: G_{k \ell} \rightarrow G_k \}_{\ell \in J_k}, k = 1, 2$.
Then $G := \pi_1(X; \ast) \cong  G_1 \ast G_2$ and thus $G_1$ and $G_2$ are finitely presentable. For $i = 1, 2$, let $(K_i, \partial K_i)$ be an Eilenberg--{Mac\,Lane} pair of type $\{ \kappa_{k \ell}: G_{k \ell} \rightarrow G_k \}_{\ell \in J_k}$ with finite $2$--skeleton. Then the pair $(K, \partial K) = (K_1 \vee K_2, \partial K_1 \cup \partial K_2)$ is an Eilenberg--{Mac\,Lane} pair of type $\{ \iota_k \circ \kappa_{k \ell}: G_{k \ell} \rightarrow G_1 \ast G_2 \}_{\ell \in J_k, k = 1, 2}$ and ${\mbox{H}}_3(K, \partial K; {\mathbb{Z}}^{\omega}) \cong {\mbox{H}}_3(K_1, \partial K_1; {\mathbb{Z}}^{\omega_1}) \oplus  {\mbox{H}}_3(K_2, \partial K_2; {\mathbb{Z}}^{\omega_2})$, where $\omega_k \in {\mbox{H}}^1(K_i; {\mathbb{Z}}^{\omega})$ is the restriction of the orientation character $\omega_X \in {\mbox{H}}^1(K; {\mathbb{Z}}^{\omega})$ for $k = 1, 2$. Hence $\mu_X = \mu_1 + \mu_2$ with $\mu_k \in {\mbox{H}}_3(K_k, \partial K_k; {\mathbb{Z}}^{\omega_k})$. It is now sufficient to show that the triple $(\{ \kappa_{k \ell}: G_{k \ell} \rightarrow G_k \}_{\ell \in J_k}, \mu_k, \omega_k)$ is realised by a ${\rm{PD}}^3$--pair $(X_k, \partial X_k)$ for $k = 1, 2$. Then the interior connected sum of $(X_1, \partial X_1)$ and $(X_2, \partial X_2)$ realises the fundamental triple of $(X, \partial X)$ and the Classification Theorem~\ref{B} implies that $(X, \partial X)$ is orientedly homotopy equivalent to $(X_1, \partial X_1) \sharp (X_2, \partial X_2)$.

For $k = 1, 2$, let $L_k$ be the functor defined by (\ref{functor}). Then
\[ C_i(K, \partial K) = L_1(C_i(K_1, \partial K_1)) \oplus L_2(C_i(K_2, \partial K_2))\]
for $i \geq 1$, and the boundary morphism $\partial_i^{K, \partial K}: C_{i+1}(K, \partial K) \rightarrow  C_{i}(K, \partial K)$ is the direct sum of $L_1(\partial_i^{K_1, \partial K_1})$ and $L_2(\partial_i^{K_2, \partial K_2})$. Thus
\[ F^2(C(K, \partial K)) = L_1(F^2(C(K_1, \partial K_1))) \oplus L_2(F^2(C(K_2, \partial K_2))) \]
and $I(G) = L_1(I(G_1)) \oplus L_2(I(G_2))$, where the canonical inclusion $L_k(I(G_k)) \rightarrowtail I(G)$ is given by $\mu \otimes {{\lambda}} \mapsto \mu {{\lambda}}$ for $\mu \in {\mathbb{Z}}[G]$, and ${{\lambda}} \in I(G_k)$ is viewed as an element of $I(G)$.

For $k = 1, 2$, let $\varphi_k: F^2(C(K_k, \partial K_k)) \rightarrow I(G_k)$ be a ${\mathbb{Z}}[G_k]$--morphism representing the class $\nu_{C(K_k, \partial K_k),2}(\mu_k)$. Then the class 
$\nu_{C(K, \partial K),2}(\mu)$ of homotopy equivalences is represented by
\[ \xymatrix{
L_1(F^2(C(K_1, \partial K_1))) \oplus L_2(F^2(C(K_2, \partial K_2))) 
\ar[d]_{L_1(\varphi_1) \oplus L_2(\varphi_2)} \ar@{=}[r] & F^2(C(K, \partial K)) \\
L_1(I(G_1)) \oplus L_2(I(G_2)) \ar@{=}[r] & I(G).} \]
It follows from the proof of the analogous proposition for the absolute case \cite{Turaev} that the ${\mathbb{Z}}[G_k]$--morphism $\varphi_k$ is a homotopy equivalence of modules for $k = 1, 2$.
Hence, by the Realisation Theorem~\ref{realthm}, $(\{ \kappa_{k \ell}: G_{k \ell} \rightarrow G_k \}_{\ell \in I_k}, \mu_k, \omega_k)$ is realised by a ${\rm{PD}}^3$--pair $(X_k, \partial X_k)$ for $k = 1, 2$.
\end{proof}

\subsection{The Boundary Connected Sum}
Let $(X_1, \partial X_1)$ and $(X_2, \partial X_2)$ be two ${\rm{PD}}^3$--pairs with connected boundary components $\{ \partial X_{1 i}\}_{i \in I_1}$ and $\{ \partial X_{2 j}\}_{j \in I_2}$ respectively. Choosing $\ell_k \in I_k, k = 1, 2$, we may assume that $\partial X_{1 \ell_1}$ and $\partial X_{2 \ell_2}$ are collared in $X_1$ and $X_2$ respectively, and that there are discs $e^2_1 \subseteq \partial X_{1 \ell_1}$ and $e^2_2 \subseteq \partial X_{2  \ell_2}$. We denote the chains corresponding to $e^2_1$ and $e^2_2$ by $e_1$ and $e_2$ respectively, and the quotient of $X_1 \coprod X_2$ obtained by identifying $e^2_1$ and $e^2_2$ via an orientation reversing map by $X_1 \coprod X_2 / \sim$. For subsets $A_i \subseteq X_i, i = 1, 2,$ we denote the image of $A_1 \coprod A_2$ under the canonical projection $\pi: X_1 \coprod X_2 \rightarrow X_1 \coprod X_2 / \sim$ by $A_1 \coprod A_2 / \sim$ . If we assume that $\partial X_{1 \ell_1}$ and $\partial X_{2 \ell_2}$ are closed surfaces, then $\partial X_{1 \ell_1} \sharp \partial X_{2 \ell_2} := (\partial X_{1 \ell_1} \setminus e^2_1) \coprod (\partial X_{2 \ell_2} \setminus e^2_2) /  \sim$ is homotopy equivalent to the connected sum of $\partial X_{1 \ell_1}$ and $\partial X_{2 \ell_2}$ as $2$--manifolds. 
\begin{lemdfn}\label{boundconsum}
Up to oriented homotopy equivalence, the pair 
\[ (X_1, \partial X_1) \natural_{\ell_1,\ell_2} (X_2, \partial X_2) := 
\big( X_1 \coprod X_2 / \sim, (\partial X_1 \setminus e_1^2) \coprod 
(\partial X_2 \setminus e_2^2) / \sim \big)\]
is a ${\rm{PD}}^3$--pair determined by $(X_i, \partial X_i), i = 1, 2,$ and $(\ell_1,\ell_2) \in I_1 \times I_2$. It is called the {\emph{$(\ell_1,\ell_2)$--boundary connected sum}} of $(X_1, \partial X_1)$ and $(X_2, \partial X_2)$.
\end{lemdfn}
\begin{proof}
Note that $X_1 \coprod X_2 / \sim$ is homotopy equivalent to the wedge $X_1 \vee X_2$. For $k = 1, 2,$ we denote the fundamental triple of $(X_k, \partial X_k)$ by $(\{ \kappa_{k \ell}: G_{k \ell} \rightarrow G_k \}_{\ell \in I_k}, \mu_k, \omega_k)$. Putting $G := \pi_1(X; \ast),$ we obtain $G = G_1 \ast G_2$. For $k = 1, 2,$ let $L_k$ be the functor defined by (\ref{functor}).

Put $(X, \partial X) := (X_1, \partial X_1) \natural_{\ell_1,\ell_2} (X_2, \partial X_2)$. Furthermore, let $D$ denote the subcomplex of $C(X)$ generated by $\pi (e_1) = \pi (e_2)$ and its boundary $\pi (\partial e_1) = \pi (\partial e_2)$, and let $\partial D$ denote the subcomplex generated by $\pi (\partial e_1) = \pi (\partial e_2)$. Then the geometric pairs $(L_1(C(X_1)), L_1(C(\partial X_1)))$, $(L_2(C(X_2)), L_2(C(\partial X_2)))$ as well as $(D, \partial D)$ are Poincar\'{e} chain pairs. By Theorem~\ref{sumthm}, 
\[(L_1(C(X_1)) + L_2(C(X_2)), L_1(C(\partial X_1) \setminus {{\Lambda}}_1[e_1]) + L_2(C(\partial X_2) \setminus {{\Lambda}}_2[e_2]))\]
is a Poincar\'{e} chain pair, where ${{\Lambda}}_i = {\mathbb{Z}}[G_i]$. Hence $(X, \partial X)$ is a ${\rm{PD}}^3$--pair with fundamental triple 
\[ (\{ \kappa: K \rightarrow G_1 \ast G_2, \ \iota_k \circ \kappa_{k \ell}: G_{k \ell} \rightarrow G_1 \ast G_2 \}_{ \ell \in I_k, \ell \neq \ell_k, k = 1, 2}, \mu_1 + \mu_2, \omega_1 + \omega_2),\]
where $K := \pi_1(\partial X_{1 \ell_1} \sharp \partial X_{2 \ell_2}; \ast)$ and $\kappa: K \rightarrow G_1 \ast G_2$ is induced by the inclusion $\partial X_{1 \ell_1} \sharp \partial X_{2 \ell_2} \rightarrow X_1 \coprod X_2 / \sim$. As the fundamental triple does not depend on the choice of discs $e^2_i, i = 1, 2$, the Classification Theorem~\ref{B} implies that, up to oriented homotopy equivalence, the ${\rm{PD}}^3$--pair $(X, \partial X)$ is uniquely determined by $(X_i, \partial X_i), i = 1, 2,$ and $(\ell_1, \ell_2) \in I_1 \times I_2$.
\end{proof}

The notion of boundary connected sums of ${\rm{PD}}^3$--pairs is consistent with that of boundary connected sums of manifolds with boundary.

Let $\{ \kappa_{k \ell}: G_{k \ell} \rightarrow G_k \}_{\ell \in J_k}$ be two $\pi_1$--systems such that $G_{k \ell}$ is a surface group for every $\ell \in J_k$ and let $(K_k, \partial K_{k})$ 
be an Eilenberg--{Mac\,Lane} pair of type $\{ \kappa_{k \ell}: G_{k \ell} \rightarrow G_k \}_{\ell \in I_k}$ 
for $k = 1, 2$. We may assume that the components of $\partial K_k$ are closed surfaces which are collared in $K_k$ and that there are discs $e^2_1 \subseteq \partial K_{1,\ell_1}$ and $e^2_2 \subseteq \partial K_{2, \ell_2}$. Let $e_1$ and $e_2$ denote the chains corresponding to $e^2_1$ and $e^2_2$ respectively and let $K_1 \coprod K_2 / \sim$ be the quotient of $K_1 \coprod K_2$ obtained by identifying $e^2_1$ and $e^2_2$ via an orientation reversing map. Put $\partial K_{1 \ell_1} \sharp \partial K_{2 \ell_2} := (\partial K_{1 \ell_1} \setminus e^2_1) \coprod (\partial K_{2 \ell_2} \setminus e^2_2) /  \sim$ and $H := \pi_1(\partial K_{1 \ell_1} \sharp \partial K_{2 \ell_2})$. Finally, let 
$\kappa: H \rightarrow G_1 \ast G_2$ be the homomorphism induced by the inclusion 
$\partial K_{1 \ell_1} \sharp \partial K_{2 \ell_2} \rightarrow K_1 \coprod K_2 / \sim.$

\begin{definition}
The free product of the $\pi_1$--systems $\{ \kappa_{k \ell}: G_{k \ell} \rightarrow G_k \}_{\ell \in J_k}$, 
$k = 1, 2$, along $G_{1 \ell_1}$ and $G_{2 \ell_2}$ is the $\pi_1$--system 
\begin{equation*}
\{ \kappa : H \rightarrow G_1 \ast G_2, \ \iota_k \circ \kappa_{k \ell}: G_{k \ell} \rightarrow G_1 \ast G_2 \}_{\ell \in I_k, \ell \neq \ell_k, k = 1, 2},
\end{equation*}
where $\iota_k: G_k \rightarrow G_1 \ast G_2, k = 1, 2$, denote the inclusions of the factors in the free product of groups. 

A $\pi_1$--system isomorphic to such a free product is said to \emph{decompose as free product along $G_{1 \ell_1}$ and $G_{2 \ell_2}$}.
\end{definition}
The proof of Lemma ~\ref{boundconsum} shows that the $\pi_1$--system of the boundary connected sum of two ${\rm{PD}}^3$--pairs $(X_k, \partial X_k)$ along $\partial X_{1 \ell_1}$ and 
$\partial X_{2 \ell_2}$ decomposes as free product along $G_{1 \ell_1}$ and $G_{2 \ell_2}$. The Classification Theorem \ref{B} and the Realisation Theorem \ref{realthm} allow us to show that the converse holds for finitely dominated pairs with non--empty aspherical boundaries in the case of $\pi_1$--injectivity, thus proving Theorem \ref{decomp2}.
\begin{proof}[Proof of Decomposition II]
Suppose the $\pi_1$--system of the finitely dominated ${\rm{PD}}^3$--pair $(X, \partial X)$ with non--empty aspherical boundary components decomposes as free product of the injective $\pi_1$--systems $\{ \kappa_{k \ell}: G_{k \ell} \rightarrow G_k \}_{\ell \in I_k}, k = 1, 2$, along $G_{1 \ell_1}$ and 
$G_{2 \ell_2}$.
Then $G := \pi_1(X; \ast) \cong  G_1 \ast G_2$ and thus $G_1$ and $G_2$ are finitely presentable. For $k = 1, 2$, let $(K_k, \partial K_{k})$ be an Eilenberg--{Mac\,Lane} pair of type $\{ \kappa_{k \ell}: G_{k \ell} \rightarrow G_k \}_{\ell \in I_k}$ with finite $2$--skeleton, let $K_1 \coprod K_2 / \sim$ be defined as above and let $L_k$ be the functor defined by (\ref{functor}), so that
\[ C_3(K, \partial K) = L_1(C_3(K_1, \partial K_1)) \oplus L_2(C_3(K_2, \partial K_2)).\]
Take a cycle $1 \otimes x$ in ${\mathbb{Z}}^{\omega} \otimes_{{\Lambda}} C_3(K, \partial K) $ representing $\mu_X$. Then there are cycles $1 \otimes x_k \in {\mathbb{Z}}^{\omega_k} \otimes_{{\Lambda}} C_3(K_k, \partial K_k)$ such that $1 {\otimes}x = 1 {\otimes}x_1 + 1 {\otimes}x_2$, where $\omega_k$ denotes the restriction of the orientation character $\omega_X \in {\mbox{H}}^1(K; {\mathbb{Z}} / 2 {\mathbb{Z}})$ to ${\mbox{H}}^1(K_k; {\mathbb{Z}} / 2 {\mathbb{Z}})$ for $k = 1, 2$. Put $\mu_k := [1 {\otimes}x_k]$ and note that
\begin{eqnarray*}
C_1(K, \partial K) & = & L_1(C_1(K_1, \partial K_1)) \oplus L_2(C_1(K_2, \partial K_2)) \\
C_2(K, \partial K) & = & L_1(C_2(K_1, \partial K_1)) \oplus L_2(C_2(K_2, \partial K_2)) \oplus {{\Lambda}} [e_1]
\end{eqnarray*}
and that the boundary morphism $\partial_1^{K, \partial K}: C_{2}(K, \partial K) \rightarrow  C_{1}(K, \partial K)$ is the direct sum of $L_1(\partial_1^{K_1, \partial K_1})$ and $L_2(\partial_1^{K_2, \partial K_2})$. Thus
\[ F^2(C(K, \partial K)) = L_1(F^2(C(K_1, \partial K_1))) \oplus L_2(F^2(C(K_2, \partial K_2))) \oplus {{\Lambda}} [e_1].\]
For $k = 1, 2$, let $\varphi_k: F^2(C(K_k, \partial K_k)) \rightarrow I(G_k)$ be a ${\mathbb{Z}}[G_k]$--morphism representing the class $\nu_{C(K_k, \partial K_k),2}(\mu_k)$.  As $\partial_1^{K, \partial K} |_{{{\Lambda}} [e_1]} = 0$, the class $\nu_{C(K, \partial K),2}(\mu)$ of homotopy equivalences yields a homotopy equivalence 
\[ \xymatrix{
L_1(F^2(C(K_1, \partial K_1))) \oplus L_2(F^2(C(K_2, \partial K_2))) 
\ar[d]_{L_1(\varphi_1) \oplus L_2(\varphi_2)} \ar@{=}[r] & F^2(C(K, \partial K)) \\
L_1(I(G_1)) \oplus L_2(I(G_2)) \ar@{=}[r] & I(G).} \]
By \cite{Turaev}, the ${\mathbb{Z}}[G_k]$--morphism $\varphi_k$ is a homotopy equivalence of modules for $k = 1, 2$ and the Realisation Theorem~\ref{realthm} implies that $(\{ \kappa_{k \ell}: G_{k \ell} \rightarrow G_k \}_{\ell \in I_k}, \mu_k, \omega_k)$ is realised by a ${\rm{PD}}^3$--pair $(X_k, \partial X_k)$ for $k = 1, 2$. The Classification Theorem~\ref{B} implies that $(X, \partial X)$ is orientedly homotopy equivalent to the $(\ell_1,\ell_2)$--boundary connected sum $(X_1, \partial X_1) \natural_{\ell_1,\ell_2} (X_2, \partial X_2)$. 
\end{proof}

\section*{Acknowledgements}
The results presented here are contained in the author's doctoral thesis for the case of finitely dominated $PD^3$--pairs. It is with gratitude that she acknowledges the support, guidance and inspiration of Jonathan Hillman, her supervisor.


\begin{thebibliography}{99}

\bibitem{Bieri-Eckmann}
R. Bieri and B. Eckmann, \emph{Relative Homology and Poincar{\'e} Duality
for Pairs}, Journal of Pure and Applied Algebra {\bf{13}} (1978), 277--319.

\bibitem{Bleile}
B. Bleile, \emph{Poincar\'e Duality Pairs of Dimension Three}, PhD Thesis, The University of Sydney, 2005

\bibitem{Bleile2}
B. Bleile, \emph{Classification of $PD^3$--pairs without finiteness conditions}, Preprint, University of New England, 2005

\bibitem{Browder}
W. Browder, \emph{Surgery on Simply Connected Manifolds}, Springer Verlag (1972).

\bibitem{Crisp2}
J. Crisp, \emph{An Algebraic Loop Theorem and the Decomposition of $\rm{PD}^3$--Pairs}, Preprint (2005), Universit\'e de Bourgogne.
 
\bibitem{Eckmann-Linnell}
B. Eckmann and P. Linnell, \emph{Poincar{\'e} Duality Groups of Dimension
Two, II}, Commentarii Mathematici  Helvetici {\bf{58}} (1983), 111--114.


\bibitem{Eckmann-Mueller}
B. Eckmann and H. M{\"u}ller \emph{Poincar{\'e} Duality Groups of Dimension
Two}, Commentarii Mathematici Helvetici {\bf{55}} (1980), 510--520.

\bibitem{Hendriks}
H. Hendriks, \emph{Obstruction Theory in 3--Dimensional Topology: An
Extension Theorem}, Journal of the London Mathematical Society (2) {\bf{16}} (1977), 160--164.

\bibitem{Hilton}
P.J. Hilton, \emph{Homotopy Theory and Duality},
Gordon and Breach (1963).

\bibitem{Johnson}
F. E. A. Johnson, \emph{Stable modules and the $D(2)$--problem}, Cambridge University PressÊ (2003)

\bibitem{KirbySiebenmann}
R.C. Kirby and L.C. Siebenmann, \emph{Foundational Essays on 
Topological Manifolds, Smoothings and Triangulations}, Annals of 
Mathematics Studies, Princeton University Press (1977).

\bibitem{Milnor1}
J.W. Milnor, \emph{Groups Which Act on $S^n$ without Fixed Points},
American Journal of Mathematics {\bf{79}} (1957), 623--630.

\bibitem{Swan}
R.G. Swan, \emph{Periodic Resolutions for Finite Groups}, 
Annals of Mathematics {\bf{72}} (1960), 267--291.

\bibitem{Turaev}
V.G. Turaev, \emph{Three Dimensional Poincar{\'e} Complexes: Homotopy
Classification and Splitting}, Math. USSR Sbornik Vol. 67 (1990),
No. 1.

\bibitem{Wall1}
C.T.C. Wall, \emph{Poincar{\'e} Complexes: I}, Annals of Mathematics 2nd Series
{\bf{86}} (1967), 213--245.

\bibitem{Wall2}
C.T.C. Wall, \emph{Finiteness Conditions for CW--Complexes}, Annals of Mathematics
2nd Series {\bf{81}} (1965), 56--69.

\bibitem{Wall3}
C.T.C. Wall, \emph{Finiteness Conditions for CW--Complexes, II}, Proceedings of the
Royal Society Series A {\bf{295}} (1966), 129--139.

\bibitem{Wall4} 
C.T.C. Wall, \emph{Surgery on Compact Manifolds}, Academic Press (1970).

\bibitem{Wall5}
C.T.C. Wall, \emph{Poincar\'e Duality in Dimension 3}, Proceedings of the Casson Fest, Geometry and Topology Monographs, Vol. 7 (2004), 1--26.

\end{thebibliography}
\end{document}